\documentclass[smallextended]{svjour3}       
\smartqed  
\usepackage{graphicx}
\usepackage{epsfig}
\usepackage{epstopdf}
\usepackage{amsfonts}
\usepackage{verbatim}
\usepackage{enumitem}
\usepackage{color}
\newenvironment{proofcite}[1]{\noindent{\bf Proof of #1.\,}}{\hfill$\Box$}
\newenvironment{proofoutlinecite}[1]{\noindent{\bf Outline for proof of #1.\,}}{\hfill$\Box$}

\newcommand{\textdef}{\emph}

\newcommand{\dotcup}{\mathbin{\mathaccent\cdot\cup}}
\newcommand{\Vexc}{The very extreme case}
\newcommand{\vexc}{the very extreme case}
\newcommand{\eps}{\varepsilon}
\newcommand{\ie}{{\em i.e.}}

\begin{document}

\title{Tiling tripartite graphs with $3$-colorable graphs: The extreme case 
       \thanks{Martin was supported in part by NSA grants H98230-05-1-0257 and H98230-08-1-0015 and by a grant from the Simons Foundation (\#353292, Ryan Martin).\\ Zhao was supported in part by NSA grants H98230-05-1-0079, H98230-07-1-0019 and NSF grants DMS-1400073 and DMS-1700622. Part of this research was done while Zhao was working at the University of Illinois at Chicago.}
}

\author{Kirsten Hogenson         \and
        Ryan R. Martin           \and
        Yi Zhao
}

\institute{K. Hogenson \at
              Colorado College, Colorado Springs, CO 80903, USA. \\
              \email{kirsten.hogenson@coloradocollege.edu}           
           \and
           R.R. Martin \at
              Iowa State University, Ames, IA 50010, USA \\
              \email{rymartin@iastate.edu}
           \and
           Y. Zhao \at
              Georgia State University, Atlanta, GA 30303, USA \\
              \email{yzhao6@gsu.edu}
}

\date{August 09, 2018}

\maketitle


\begin{abstract}
   There is a sufficiently large $N\in h\mathbb{N}$
   such that the following holds. If $G$ is a tripartite graph with $N$ vertices in each vertex class such
   that every vertex is adjacent to at least  $2N/3+2h-1$ vertices
   in each of the other classes, then $G$ can be tiled perfectly by copies
   of $K_{h,h,h}$.  This extends work by two of the authors [Electron. J. Combin, 16(1), 2009] and also
   gives a sufficient condition for tiling by any fixed
   3-colorable graph. Furthermore, we show that $2N/3+2h-1$ in our result can not be replaced by $2N/3+ h-2$
   and that if $N$ is divisible by $6h$, then we can replace it with the value $2N/3+h-1$ and this is tight.

\end{abstract}


\section{Introduction}
Let $H$ be a graph on $h$  vertices, and let $G$ be a graph on $n$ vertices.  An $H$-tiling of $G$ is a subgraph of $G$ which consists of vertex-disjoint copies of $H$ and a {\em perfect $H$-tiling}, or {\em $H$-factor}, of $G$ is an $H$-tiling consisting of $\lfloor n/h\rfloor$ copies of $H$.
The celebrated Hajnal-Szemer\'edi Theorem \cite{HaSz} says that each $n$-vertex graph $G$ with $\delta(G)\ge (r-1)n/r$ contains a $K_r$-factor. 
(Corr\'adi and Hajnal~\cite{CoHa} proved the case $r=3$.)  Using Szemer\'edi's regularity lemma \cite{Sz}, 
Alon and Yuster \cite{AlYu1,AlYu2} obtained results on $H$-tiling for arbitrary $H$. Their results were improved substantially \cite{Komlos,KSSz-AY,KuhnOsthus,AliYi}, in particular, K\"uhn and Osthus~\cite{KuhnOsthus} determined the minimum degree threshold for $H$-factors for arbitrary $H$ up to an additive constant,
see the survey 
\cite{KuOs-survey} for details.

In this paper, we consider multipartite tiling, which restricts $G$
to be an $r$-partite graph. For $r=2$, this is an immediate consequence of the K\"{o}nig-Hall Theorem (e.g. see \cite{Bollobas}). Wang \cite{Wang98} considered
$K_{s,s}$-factors in bipartite graphs for all $s>1$; Zhao \cite{Zhao} gave the best possible minimum degree condition for this problem. With the exception of one case, Hladk\'y and Schacht \cite{HlSc} found best possible minimum degree conditions for $K_{s,t}$-factors in bipartite graphs with $s<t$; the last case was settled by Czygrinow and DeBiasio \cite{CzDeBi}.  Later, Bush and Zhao \cite{BuZh} considered tiling bipartite graphs with an arbitrary graph $H$.

For a tripartite graph $G=(A,B,C;E)$, the graphs induced by $(A,B)$, $(A,C)$ and $(B,C)$ are called the \textdef{natural bipartite subgraphs} of $G$.  Let ${\cal G}_r(N)$ be the family of $r$-partite graphs with $N$ vertices in each partition set. Such a graph is called \textdef{balanced} because the number of vertices in each partition set is the same. In an $r$-partite graph $G$, $\delta^*(G)$ stands\footnote{In \cite{MZ}, $\bar{\delta}$ was used in place of $\delta^*$.} for the minimum degree over all natural bipartite subgraphs of $G$.

There are two classes of multipartite graphs that we will reference in this paper.  One is $\Gamma_k$, which is in ${\cal G}_k(k)$.  The vertices of $\Gamma_k$ are $h^{(j)}_i$, $i=1,\ldots,k$ and $j=1,\ldots,k$, and the adjacency rules are as follows: $h^{(j)}_i\sim h^{(j')}_{i'}$ iff $i\ne i'$, $j\ne j'$, and either $j$ or $j'$ is in $\{1,\dots,k-2\}$.  Also, $h^{(k-1)}_i\sim h^{(k-1)}_{i'}$ and $h_i^{(k)}\sim h_{i'}^{(k)}$ for $i\ne i'$.  The other graph is $\Theta_{r\times n}$, which is in ${\cal G}_r(n)$.  The vertices of $\Theta_{r\times n}$ are $a^{(j)}_{i}$, $i=1,\ldots,r$ and $j=1,\ldots,n$ such that $a^{(j)}_{i}\sim a^{(j')}_{i'}$ if and only if $i\neq i'$ and $j\neq j'$.  We will also discuss the so-called blow-ups of these graphs.  The \emph{blow-up graph}, $G(N)$, for a graph $G$ is obtained by replacing each edge of $G$ with a copy of $K_{N,N}$ and replacing each non-edge by an $N\times N$ bipartite graph with no edges.

\begin{figure}\label{fig:GAMMA}
   \begin{center}
      \includegraphics[width=2in]{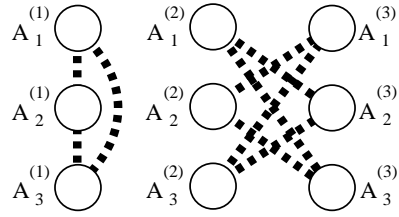}
   \end{center}
   \caption{The diagram for the blow-up of $\Gamma_3$, with the vertex classes in rows and the dotted lines representing \emph{nonedges}. Note that $\delta^*(\Gamma_3)=2$ but $\Gamma_3$ has no $K_3$-factor.}
\end{figure}

In addition to the bipartite results discussed above, there have also been a number of results on multipartite graphs with $r\geq3$, many of which were inspired by a conjecture of Fischer.
Fischer \cite{Fischer} conjectured that if $G\in{\cal G}_r(N)$ satisfies $\delta^*(G) \geq \frac{r-1}{r}N$, then $G$ contains a $K_r$-factor. However, if $r$ and $N/r$ are odd integers, then Catlin \cite{Catlin} had earlier given an example of a graph without a $K_r$-factor where $\delta^*(G) = \frac{r-1}{r}N$.
In \cite{MM}, Magyar and Martin proved that, for large $N$, this graph is a unique counterexample to Fisher's conjecture
for $r=3$ by showing that if $N$ is a sufficiently large odd multiple of
3, the blow-up graph $\Gamma_3(N/3)\in {\cal G}_3(N)$ (see Figure~\ref{fig:GAMMA}) is the unique graph with $\delta^*(G)\geq 2N/3$ and no $K_3$-factor. The conjecture of Fischer can be modified to exclude this case. This gives the following Corr\'adi-Hajnal-type result.

\begin{theorem}[\cite{MM}]
Let $G\in {\cal G}_3(N)$ have $\delta^*(G)\ge (2/3)N$.  If $N\geq N_0$ for some absolute constant $N_0$, then $G$ contains a $K_3$-factor or $G=\Gamma_3(N/3)$ for $N/3$ an odd integer.

\label{thm:MM}
\end{theorem}

Martin and Szemer\'edi \cite{MSz} proved a quadripartite version of the Hajnal-Szemer\'edi Theorem. Han and Zhao \cite{HanZh} reproved the results of \cite{MM,MSz} by using the absorbing method.  
An approximate version of the multipartite Hajnal-Szemer\'edi Theorem was given by Csaba and Mydlarz \cite{Csaba}.
Keevash and Mycroft \cite{KeMy1} and independently Lo and Markstrom \cite{LoMa} confirmed Fischer's conjecture asymptotically, and finally Keevash and Mycroft \cite{KeMy} proved the modified Fischer conjecture exactly for any sufficiently large graph.
More recently, an asymptotic multipartite version of the the Alon-Yuster Theorem was proved by Martin and Skokan \cite{MSk}.


It was shown \cite[Theorem 1.2]{MZ} that in the tripartite case, $2/3$ is the correct coefficient of $N$ required to have a $K_{h,h,h}$-factor.
\begin{theorem}[\cite{MZ}]\label{thm:MZ}
   For any positive real number $\gamma$ and any positive integer $h$, there is $N_0$ such that the following holds.  Given an integer $N\geq N_0$ such that $N$ is divisible by $h$, if $G$ is a tripartite graph with $N$ vertices in each vertex class such that every vertex is adjacent to at least $(2/3+\gamma)N$ vertices in each of the other classes, then $G$ contains a $K_{h,h,h}$-factor.
\end{theorem}

Let $f(N,h)$ be the smallest integer $f$ such that every balanced tripartite graph $G$ on $3N$ vertices with $\delta^*(G)\geq f$
contains a $K_{h,h,h}$-factor.  Our main result is the following more precise theorem.

\begin{theorem}
Fix a positive integer $h$ and let $N$ be sufficiently large.
If $h\geq 2$ and $N=(6q+r)h$ with $0\leq r<6$, then
$$ \begin{array}{lcccll}
   & & f(N,h) & = & \frac{2N}{3}+h-1, & \mbox{ if $r=0$;} \\
   h\left\lceil\frac{2N}{3h}\right\rceil+h-2
   & \leq & f(N,h) & \leq &
   h\left\lceil\frac{2N}{3h}\right\rceil+h-1,
   & \mbox{ if $r=1,2,4,5$;} \\
   \frac{2N}{3}+h-1 & \leq & f(N,h) & \leq & \frac{2N}{3}+2h-1, &
   \mbox{ if $r=3$.}
   \end{array} $$
 \label{thm:trithm}
\end{theorem}

So, the result is tight when $6h$ divides $N$, almost tight unless $N/h$ is an odd multiple of $3$ and, in the worst case, the upper and lower bounds
differ by $h$. We are not sure whether the upper or lower bounds of Theorem~\ref{thm:trithm} are correct in the cases when they are not equal.

Clearly the complete tripartite graph $K_{h,h,h}$ can itself be perfectly
tiled by any 3-colorable graph on $h$ vertices. Since $f(N,h)\leq\frac{2N}{3}+2h-1$ whenever $N$ is divisible by $h$, we have the
following corollary.

\begin{corollary}\label{cor:main}
   Let $H$ be a $3$-colorable graph of order $h$. There
exists a positive integer $N_0$ such that if $N\geq N_0$ and $N$
divisible by $h$, then every $G\in {\cal G}_3(N)$ with
$\delta^*(G)\ge \frac{2N}{3} + 2h -1$ contains an $H$-factor.
\end{corollary}

The lower bound for $f(N,h)$ in Theorem~\ref{thm:trithm} is due to two constructions, one which is from~\cite{MZ} and another which is similar.  They are stated in Proposition~\ref{prop:lb}, and proven in Section~\ref{sec:lb}.
\begin{proposition}
Fix a positive integer $h\geq 2$.  There exists an $N_0$ such that
\begin{enumerate}[label={(\arabic*)}]
   \item if $N\geq N_0$, $h\mid N$ and $N/h$ is divisible by $3$, then there is a graph $G_2\in {\cal G}_3(N)$ with no $K_{h,h,h}$-factor and $\delta^*(G_2)\geq 2N/3+h-2$; and \label{it:Theta3aExtreme}
   \item if $N\geq N_0$, $h\mid N$ and $N/h$ is not divisible by $3$, then there is a graph $G_3\in {\cal G}_3(N)$ with no $K_{h,h,h}$-factor and $\delta^*(G_3)\geq h\left\lceil\frac{2N}{3h}\right\rceil+h-3$. \label{it:Theta3bExtreme}
\end{enumerate}
\label{prop:lb}
\end{proposition}

As to the upper bound, we use Theorem~\ref{thm:trifuzz} (Theorem 1.4 from~\cite{MZ}) to take care of the main case.
For vertex sets $A$ and $B$, let $d(A,B) :=\frac{e(A, B)}{|A||B|}$ denote the \emph{density} of $A$ and $B$.  Before we can prove the main case, we need the following definition.
\begin{definition}
Given $\alpha>0$, we say that
$G=\left(V_1,V_2,V_3;E\right)\in{\cal G}_3(N)$ is
\textbf{$\alpha$-extremal} when there are three sets
$A_1, A_2, A_3$ such that $A_i\subseteq V_i$, $|A_i|=\lfloor N/3\rfloor$ for all $i$ and
$d(A_i,A_j)\le \alpha$
for $i\neq j$.
\end{definition}
If $G$ is $\alpha$-extremal and $\delta^*(G)\geq 2N/3$, then for $i\neq j$, the pair $(A_i,V_j- A_j)$ is a very dense bipartite graph.  Thus, we expect most members of our $K_{h,h,h}$-factor with vertices in $A_i$ to have $h$ vertices in $A_i$ and the remaining $2h$ vertices in $\left(V_j- A_j\right)\cup \left(V_k- A_k\right)$, where $\{j,k\}=[3]-\{i\}$.

\begin{theorem}[\cite{MZ}]
   Given any positive integer $h$ and any $\alpha>0$, there exists an $\eps>0$ and an integer
   $N_0$ such that whenever $N\geq N_0$, and $h$
   divides $N$, the following occurs: If $G\in{\cal G}_3(N)$ satisfies $\delta^*(G)\geq
   (2/3-\eps)N$, then either $G$ contains a $K_{h,h,h}$-factor or $G$ is
   $\alpha$-extremal. \label{thm:trifuzz}
\end{theorem}

Hence, for the upper bound, it suffices to assume that $G\in{\cal G}_3(N)$ is $\alpha$-extremal.  The proof, given in Section~\ref{sec:extreme}, is detailed and involves a case analysis.  Moreover, it requires the definition of a particular structure we call \emph{\vexc}, which we deal with in Section~\ref{sec:vexc}. This definition is given below, but roughly, it means that the graph looks like $\Gamma_3(N)$.

\begin{definition}\label{def:vexc}
A balanced tripartite graph $G$ on $3N$ vertices is in \textbf{\vexc}~if the
following occurs:  First, there are integers $h,q$ such that
$N=(6q+3)h$.  Second, there are sets $U_i^{(j)}\subseteq V_i$
for $i,j\in\{1,2,3\}$, each with size at least $2qh+1$, such that if
$v\in U_i^{(j)}$ then $v$ is nonadjacent to at most $3h-3$
vertices in $U_{i'}^{(j')}$ whenever $(h_i^{(j)},h_{i'}^{(j')})$ is an edge in the graph
$\Gamma_3$.
\end{definition}

Note that we use different language for $\alpha$-extremal and \vexc~because the definition of $\alpha$-extremal requires a parameter, whereas \vexc~does not.

Now that we have defined \vexc, we can formally state the upper bound theorem as follows:
\begin{theorem}\label{thm:extreme}
   Fix $h\geq 2$.  Let $N\in h\mathbb{N}$ be sufficiently large and assume $G\in{\cal G}_{3}(N)$.
   If $\delta^*(G)\geq h\left\lceil\frac{2N}{3h}\right\rceil+h-1$, then $G$ has a $K_{h,h,h}$-factor or $G$ is in~\vexc.
   If $G$ is in~\vexc~and $\delta^*(G)\geq h\left\lceil\frac{2N}{3h}\right\rceil+2h-1$, then $G$ has a $K_{h,h,h}$-factor.
\end{theorem}


\section{Lower bound}
\label{sec:lb}

First, we need a lemma (Lemma 2.1 in~\cite{MZ}) which permits sparse tripartite graphs with no triangles and with no quadrilaterals in its natural bipartite subgraphs:
\begin{lemma}\label{lem:Sidon}
For each integer $d\geq 0$, there exists an $n_0$ such that, if $n\geq n_0$, there exists a balanced tripartite graph, $Q(n,d)$ on $3n$ vertices such that each of the $3$ natural bipartite subgraphs are $d$-regular with no $C_4$ and $Q(n,d)$ has no $K_3$.
\end{lemma}

Finally, we prove the lower bound given in Proposition~\ref{prop:lb}. Note that Proposition~\ref{prop:lb}\ref{it:Theta3aExtreme} is proved by Proposition 1.5 in~\cite{MZ}, so here we only address Proposition~\ref{prop:lb}\ref{it:Theta3bExtreme}.\\

\begin{proofcite}{Proposition~\ref{prop:lb}\ref{it:Theta3bExtreme}} \\
%
%
Let $h\geq 3$ and $N=(3q+r)h$ so that, in this case, $r\in\{1,2\}$.  Let $G_3$ be defined such that $V_i=A_i^{(1)}\dotcup A_i^{(2)}\dotcup A_i^{(3)}$ (the notation $\dotcup$ emphasizes that it is a disjoint union of sets) in which \textdef{column $j$} is defined to be the triple of the form $(A_1^{(j)},A_2^{(j)},A_3^{(j)})$.  Let the graph in column 1 be
$Q(qh+rh-1,rh+h-4)$ where $rh+h-4\geq 2$, the graph in column 2 be $Q(qh,h-3)$ and the graph in column 3 be
$Q(qh+1,h-2)$.  If two vertices are in different columns and different vertex-classes, then they are adjacent.  It is easy to verify that $\delta^*(G_3)=2qh+rh+(h-3)=h\lceil (2N)/(3h)\rceil+h-3$.  Suppose, by way of contradiction, that $G_3$ has a $K_{h,h,h}$-factor.

If a copy of $K_{h,h,h}$ has vertex classes $U_1,U_2,U_3$, then $U_i\subseteq V_j$ for some $j$.  Since there are no triangles in any column and no $C_4$'s in the natural bipartite subgraphs of a column, the intersection of a copy of $K_{h,h,h}$ with a column is either a star with all leaves in the same vertex-class, or a set of vertices in the same vertex-class.  So each copy of $K_{h,h,h}$ has at most $h+1$ vertices in column 1 and at most $h$ vertices in each of column 2 and column 3.

There are three cases for a copy of $K_{h,h,h}$.  Case 1 has $h$ vertices in each column.  Case 2 has $h+1$ vertices in column 1, $h-1$ vertices in column 2 and $h$ vertices in column 3.  Case 3 has $h+1$ vertices in column 1, $h$ vertices in column 2 and $h-1$ vertices in column 3.

In Cases 1 and 2, since $G_3$ contains no $K_{1,h-1}$ in column 3, having $h$ vertices of a $K_{h,h,h}$ in column $3$ implies that all of them are in the same vertex class.  In Case 3, since $G_3$ has no $K_{1,h-1}$ in column 2, having $h$ vertices in column 2 means that all are in the same vertex-class.  Since $h+1$ vertices in column 1 means that they form a star, the remaining $h-1$ vertices in column 3 must be in the same vertex-class (the same vertex-class as the center of the star).  Hence, the intersection of any copy of $K_{h,h,h}$ with column 3 is contained within a single vertex-class.  Therefore, the number of copies of $K_{h,h,h}$ in the $K_{h,h,h}$-factor of $G_3$ is at least $3\left\lceil\frac{qh+1}{h}\right\rceil=3q+3$, a contradiction because the factor has exactly $3q+r\leq 3q+2$ copies of $K_{h,h,h}$.

Next consider the case when $h=2$ and $N=2(3q+r)$ with $r\in\{1,2\}$. Let $G_3$ be defined such that the graph in column 1 is $Q(2q+1,0)$, but all other possible edges in $G_3$ are present. It is easy to verify that $\delta^*(G_3)=4q+2r-1 = 2\lceil N/3\rceil -1$.  Suppose, by way of contradiction, that $G_3$ has a $K_{2,2,2}$-factor.  The intersection of one copy of $K_{2,2,2}$ with column 1 must be contained within a single vertex class and can contain at most $2$ vertices.  So at least $3\left\lceil\frac{2q+1}{2}\right\rceil = 3q+3$ copies of $K_{2,2,2}$ are needed to cover all of column 1.  This is a contradiction, because the factor has exactly $3q+r\leq 3q+2$ copies of $K_{2,2,2}$.
%
%
\end{proofcite}


\section{The extreme case}
\label{sec:extreme}

Throughout Section~\ref{sec:extreme}, assume that $G$ is
minimal, \ie, no edge of $G$ can be deleted so that the minimum
degree condition still holds. As we complete the proof of Theorem~\ref{thm:trithm} by proving Theorem~\ref{thm:extreme},
we will develop the usual hierarchy of
constants:
$$ \alpha\ll\alpha_1\ll\alpha_2\ll
   \alpha_3\ll\alpha_4\ll\alpha_5\ll 1-\theta\ll h^{-1} . $$~\\

\begin{figure}
   \begin{center}
      \includegraphics[width=2in]{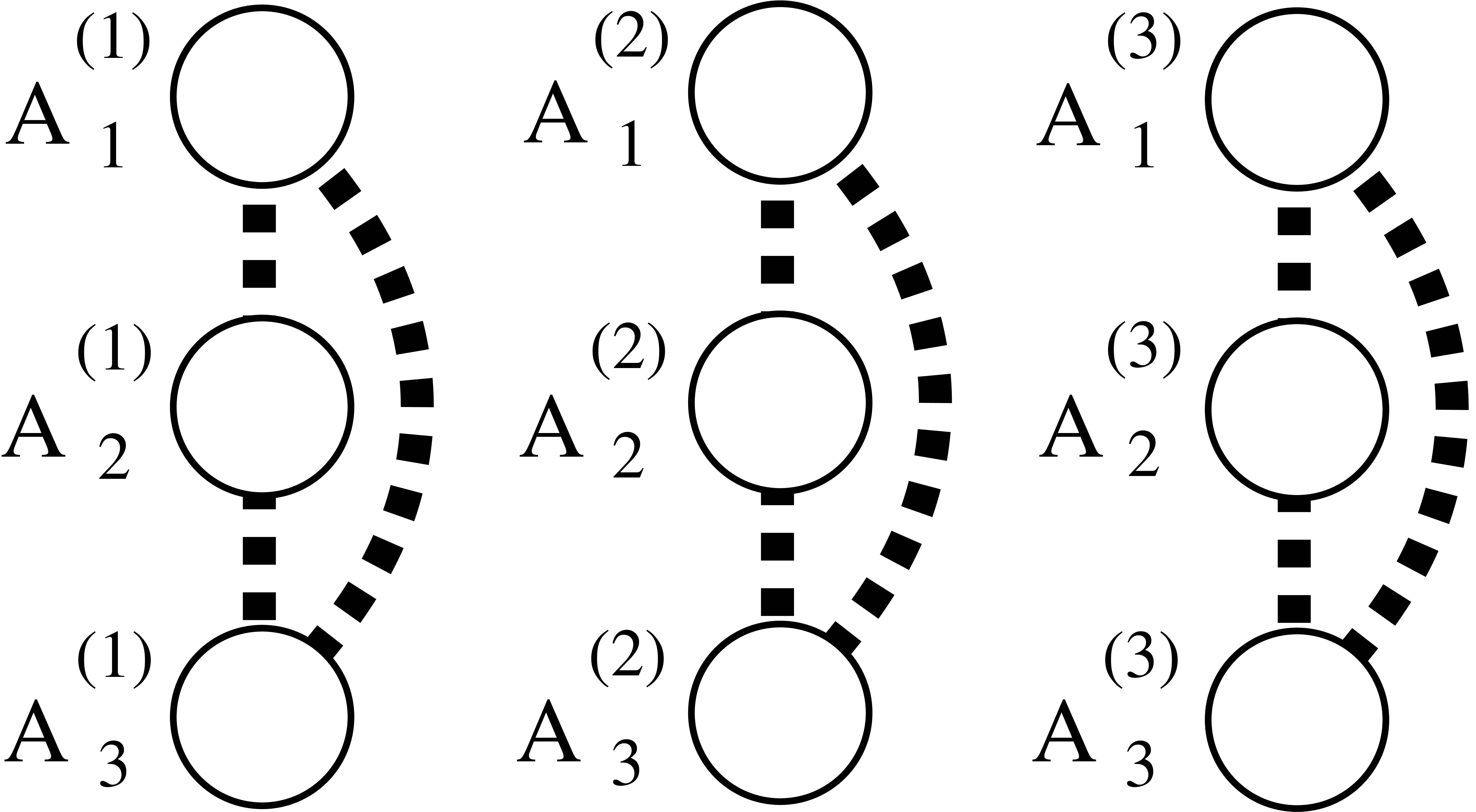}
   \end{center}
   \caption{The diagram for $\Theta_{3\times 3}$, with the vertex classes in rows and the dotted lines representing \emph{nonedges}.}
   \label{fig:THETA}
\end{figure}

\noindent\textbf{Brief outline of the proof.} There are 4 parts to the proof.  Part 1 begins with $G$ being $\alpha$-extremal and seeks a $K_{h,h,h}$-factor.  If such a tiling is not found in $G$, we deduce that $G$ looks like the graph in Figure~\ref{fig:figTWO} and move to Part 2. We again seek a $K_{h,h,h}$-factor in $G$, and if it is not found, then we move on to Part 3 which addresses the two potential structures $G$ must have.  In Part 3a, $G$ is approximately $\Theta_{3\times 3}(N/3)$.  (See Figure~\ref{fig:THETA}.)
In Part 3b, $G$ is approximately $\Gamma_3(N/3)$.  (See Figure~\ref{fig:GAMMA}.) Proofs for the lemmas and propositions stated in this section are deferred until Section~\ref{sec:lemmas}. ~\\

The following definition will come into play as we describe the structure of $G$.
\begin{definition}
   For $\delta$, $0<\delta<1$, a graph $H$ and positive integer
   $m$, we say a graph $G$ is {\em $\delta$-approximately $H(m)$}
   if $V(G)$ can be partitioned into $|V(H)|$ nearly-equally sized
   pieces, each of size $m$ or $m+1$, corresponding to a vertex of $H$ so that for
   vertices $v,w\in V(H)$ with $v\not\sim_H w$, the parts of
   $V(G)$ corresponding to $v$ and $w$ have pairwise density
   less than $\delta$.
\end{definition}

Note that if $v\sim_H w$, we do not require that the parts of $V(G)$ corresponding to $v$ and $w$ have pairwise density close to 1.

We will assume for Parts 1, 2, 3a and 3b
(Sections~\ref{sec:appfigONE},~\ref{sec:appfigTWO},~\ref{sec:theta} and~\ref{sec:gamma}, respectively) that $\delta^*(G)\geq
h\left\lceil\frac{2N}{3h}\right\rceil+h-1$.  This takes care of everything except for \vexc, which we will consider in Section~\ref{sec:vexc}.  For this last part, we will require
$\delta^*(G)\geq h\left\lceil\frac{2N}{3h}\right\rceil+2h-1$ to complete the proof.

\subsection{Part~1: The basic extreme case}
\label{sec:appfigONE}

For Part~1, we will prove that either a $K_{h,h,h}$-factor exists in
$G$, or $G$ is in Part~2.

Let $A_i\subset V_i$ for $i=1,2,3$ be the three pairwise
sparse sets given by the definition of $\alpha$-extremal and
$B_i=V_i- A_i$ for $i=1,2,3$. Recall that $|A_i| = \lfloor N/3\rfloor$, so $|B_i|= \lceil 2N/3\rceil$.
We then define
$\widetilde{A}_i$ to be the set of \emph{typical} vertices with respect to
$A_i$, $\widetilde{B}_i$ to be the set of \emph{typical} vertices with respect to
$B_i$, and $\widetilde{C}_i$ are what remain.  Formally, for $i=1,2,3$,
\begin{eqnarray*}
   \widetilde{A}_i & = & \left\{x\in V_i : \forall j\neq i, ~\deg_{A_j}(x)\leq\alpha_1|A_j|
   \right\} \\
   \widetilde{B}_i & = & \left\{y\in V_i : \forall j\neq i,
   ~\deg_{A_j}(y)\geq(1-\alpha_1)|A_j| \right\} \\
   \widetilde{C}_i & = & V_i- \left(\widetilde{A}_i\cup \widetilde{B}_i\right).
\end{eqnarray*}

Let $\{i,j,k\}=\{1,2,3\}$. Using these definitions, the fact that $G$ is $\alpha$-extremal and the bound on $\delta^*$, and the fact that every member of $A_i -\widetilde{A}_i$ is adjacent to at least an $\alpha_1$ proportion of either $A_j$ or $A_k$, we obtain the following:
\begin{eqnarray*}
    |A_i - \widetilde{A}_i| \cdot \alpha_1 \lfloor N/3\rfloor &<& e(A_i-\widetilde{A}_i, A_j) + e(A_i-\widetilde{A}_i, A_k) \leq e(A_i, A_j) + e(A_i, A_k) \\
    &\leq & 2\alpha \lfloor N/3\rfloor^2,
\end{eqnarray*}
and
\begin{eqnarray*}
    \left\lceil \frac{2N}{3} \right\rceil |A_j| &\leq& e(V_i,A_j) \\
        &\leq& \alpha|A_i||A_j| + |B_i-\widetilde{B}_i|(1-\alpha_1)|A_j| + |B_i\cap \widetilde{B}_i||A_j| \\
        &=& \alpha|A_i||A_j| + |B_i||A_j| -\alpha_1 |B_i-\widetilde{B}_i||A_j|.
\end{eqnarray*}
As a result, we have that $|A_i- \widetilde{A}_i|\leq 2(\alpha/\alpha_1)\lfloor N/3 \rfloor$ and $|B_i-
\widetilde{B}_i|\leq 2(\alpha/\alpha_1)\lfloor N/3 \rfloor$.  So, with $\alpha_1=\alpha^{1/3}$ and $\alpha_2=4\alpha_1^2$, we get the following bounds for $|\widetilde{A}_i|$ and $|\widetilde{B}_i|$:
$$ (1-\alpha_2)\lfloor N/3 \rfloor \leq |\widetilde{A}_i| \leq (1+\alpha_2)\lfloor N/3 \rfloor $$
and
$$ (1-\alpha_2)\lceil 2N/3 \rceil \leq |\widetilde{B}_i| \leq (1+\alpha_2)\lceil 2N/3 \rceil. $$

\paragraph{Step~1: Adjusting the sizes of the $\widetilde{A}_i$ sets.}

Let $N=(3q+r)h$ with $r\in\{0,1,2\}$ and $T=h\left\lfloor N/(3h)\right\rfloor$.

Without loss of generality, assume that $|\widetilde{A}_1| \geq |\widetilde{A}_2| \geq |\widetilde{A}_3|$.  For $i=1,2,3$, define $a_i = T+h$ if $i\leq r$; otherwise, $a_i=T$.  If $|\widetilde{A}_i| >a_i$, then we will move $|\widetilde{A}_i|-a_i$ vertices of $\widetilde{A}_i$ to $\widetilde{B}_i$ by applying
Lemma~\ref{lem:stars} below, which is proved in Section~\ref{sec:lemmas}.  It is applied several times throughout this paper to different sets.

\renewcommand{\labelenumi}{(\arabic{enumi})}
\begin{lemma}
   Let us be given $\epsilon_{\ref{lem:stars}}>0$ and a positive integer $M$.
   \begin{enumerate}
      \item Let $({A}_1,{A}_2;E)$ be a bipartite graph such that
      every vertex in ${A}_2$ is adjacent to at least $d_1$
      vertices in ${A}_1$.  Suppose further that
      $\left||{A}_i|-M\right|<\epsilon_{\ref{lem:stars}} M$ and $d_i<\epsilon_{\ref{lem:stars}} M$ for
      $i=1,2$.

      Provided $\epsilon_{\ref{lem:stars}}<\left((h+1)h\right)^{-1}$, there is a family of $\max\{0,d_1-h+1\}$ vertex-disjoint copies of $K_{1,h}$ all of whose centers lie in ${A}_1$.
      \label{lem:stars:bi}

      \item Let $({A}_1,{A}_2,{A}_3;E)$ be a tripartite graph such that
      every vertex not in ${A}_i$ is adjacent to at least $d_i$
      vertices in ${A}_i$, for $i=1,2,3$.  Suppose further that
      $\left||{A}_i|-M\right|<\epsilon_{\ref{lem:stars}} M$ and $d_i<\epsilon_{\ref{lem:stars}} M$ for
      $i=1,2,3$. \label{lem:stars:tri}

      Provided $\epsilon_{\ref{lem:stars}}<(2(h+2)(h+1)h)^{-1}$, there
      is a family of $\max\{0,d_i-h+1\}$ vertex-disjoint copies of $K_{1,h}$ all of whose centers lie in ${A}_i$ and leaves lie in ${A}_{i+1}$ (index arithmetic is modulo 3).
   \end{enumerate}
   \label{lem:stars}
\end{lemma}

Since $\left\lceil\frac{2N}{3h}\right\rceil + \left\lfloor\frac{N}{3h}\right\rfloor = \frac{N}{h}$, we have $h\left\lceil\frac{2N}{3h}\right\rceil +T = N$.
As $\delta^*(G)\geq
h\left\lceil\frac{2N}{3h}\right\rceil+h-1 \geq N-T+h-1$, we can guarantee that each vertex not in
$V_i$ is adjacent to at least $|\widetilde{A}_i|-T+h-1$
vertices in $\widetilde{A}_i$. So we apply
Lemma~\ref{lem:stars}{\it (\ref{lem:stars:tri})} to the graph induced by $(\widetilde{A}_1,\widetilde{A}_2,\widetilde{A}_3)$, with $d_i=
|\widetilde{A}_i|-T+h-1$, $\epsilon_{\ref{lem:stars}}=\alpha_2$, and $M = N/3$. This will construct stars with the property that there are exactly enough centers in $\widetilde{A}_i$ such that, when removed, the resulting set has its size bounded above by $a_i$, which is either $T$ or $T+h$, depending on the case. Let $Z_i$ denote the set of these centers and move the desired number of vertices of $Z_i$ from $\widetilde{A}_i$ into $\widetilde{B}_i$.


If $|\widetilde{A}_i|<a_i$, then we will move $a_i-|\widetilde{A}_i|$ vertices of $\widetilde{B}_i \cup \widetilde{C}_i$ to $\widetilde{A}_i$, as follows.

For a subgraph $K_{1,h,h}$, with $h\geq 2$, define the {\em
center} to be the vertex that is adjacent to all others. We will refer to the remaining vertices as {\em leaves}, although
their degree is $h+1$.

In $B := \bigcup_{i=1}^3 \left(\widetilde{B}_i \cup \widetilde{C}_i\right)$, we will find vertex-disjoint copies of $K_{1,h,h}$ such that each of
$\max\{a_i-|\widetilde{A}_i|,0\}$ copies has its center
vertex in $\widetilde{B}_i\cup
\widetilde{C}_i$ for $i\leq r$
and such that each of $a_i-|\widetilde{A}_i|$ copies has its center vertex in $\widetilde{B}_i\cup
\widetilde{C}_i$ otherwise. This will be accomplished with Lemma~\ref{lem:superstars}, which is proved in Section~\ref{sec:lemmas}. It is applied several times throughout this paper to slight variations of the sets $\widetilde{B}_i$.

\begin{lemma}
   Given $\delta>0$, there exists an
   $\epsilon_{\ref{lem:superstars}}=\epsilon_{\ref{lem:superstars}}(\delta)>0$ such that the
   following occurs:

   Let $({B}_1,{B}_2,{B}_3;E)$ be a tripartite graph on $6M$ vertices such that for all
   $i\neq j$, each vertex in ${B}_i$ is adjacent to at least $(1-\epsilon_{\ref{lem:superstars}})M$ vertices in ${B}_j$.  Furthermore, $\left||{B}_i|-2M\right|<\epsilon_{\ref{lem:superstars}} M$.

   If $({B}_1,{B}_2,{B}_3;E)$ contains no copy of $K_{1,h,h}$ with 1 vertex in ${B}_1$, and $h$ vertices in each of ${B}_2$ and ${B}_3$,
   then the graph $({B}_1,{B}_2,{B}_3;E)$ is $\delta$-approximately $\Theta_{3\times 2}(M)$.
   \label{lem:superstars}
\end{lemma}

Lemma~\ref{lem:superstars} can be repeatedly applied to $B$ at most
$\lceil\alpha_2(N/3)\rceil$ times with $\delta=\alpha_3$, $\alpha_2\ll\epsilon_{\ref{lem:superstars}}$ and $M=T$.  Each time, either a $K_{1,h,h}$ is found and removed, or the current incarnation of $B$ is $\alpha_3$-approximately $\Theta_{3\times 2}(T)$ and we stop applying the lemma.  When we are finished applying Lemma~\ref{lem:superstars}, add the center vertices of the $K_{1,h,h}$ subgraphs to the appropriate sets
$\widetilde{A}_i$.  Put the leaves back into $B$ and denote the result as $B = (B_1,B_2,B_3;E)$.

If necessary, place vertices from $\widetilde{C}_i$ into the set $\widetilde{A}_i$, for $i=1,2,3$, so that the resulting set, relabeled as $A^{(1)}_{i}$, is of size $a_i$ and $\sum_{i=1}^3|A^{(1)}_{i}|=N$.

\paragraph{Step~2: Finding a $K_{h,h}$-factor in $B$.}

Now we try to find a $K_{h,h}$-factor among the remaining vertices in $B$ with the goal of extending each $K_{h,h}$ into a $K_{h,h,h}$ using vertices in $A_1^{(1)}\cup A_2^{(1)}\cup A_3^{(1)}$.  Before we do so, however, we need to address the following concerns:

\begin{itemize}
   \item Vertices in copies of $K_{1,h,h}$ where the center vertex is in some $A^{(1)}_{i}$ must be in a specified copy of $K_{h,h}$ in $B$.
   \item Recall that $Z_i$ is the set of centers of $h$-stars which were found in Step~1.  If $v\in Z_i$ is the center of a $K_{1,h}$ with leaves in $A^{(1)}_{k}$, then $v$ will be \textdef{assigned to} $B_j$, where $\{j\}=\{1,2,3\}-\{i,k\}$.  This means that $v$ will be adjacent to vertices in $B_j$ in a $K_{h,h}$ in $B$.
   \item For $\{i,j,k\}=\{1,2,3\}$, vertices $v\in \widetilde{C}_i$ will be \textdef{assigned to} $B_j$ or $B_k$, respectively.  This means that $v$ will be adjacent to either $h$ vertices in $B_j$ or $h$ vertices in $B_k$ in a $K_{h,h}$ to be formed in $B$ together with $h-1$ vertices in $B_i$.  We know this can be accomplished because if $v\in \widetilde{C}_i$, then we may assume, without loss of generality, that $v$ is adjacent to at least $\alpha_1|A_j|$ vertices in $A_j$.
\end{itemize}

Moreover, because all but a $\alpha_2$-proportion of the sets $A_i$ and $B_i$ are typical, we have that $|\widetilde{C}_i|\leq \alpha_2|A_i|+\alpha_2|B_i|\leq 3\alpha_2T$.  Recall that we applied Lemma~\ref{lem:stars} with $d_i=|\widetilde{A}_i|-T+h-1$.  Thus $|Z_i|\leq\alpha_2|A_i|+h-1\leq 2\alpha_2T$ and there are at most $\alpha_2|A_i|+h\leq 2\alpha_2T$ copies of $K_{1,h,h}$ with the center vertex in a given $A^{(1)}_{i}$.

Lemma~\ref{lem:partition} is proved in Section~\ref{sec:lemmas}. We will apply it to an adjusted $B$ where we know from Step~1 there are copies of $K_{h,h}$ which must belong to any $K_{h,h}$-factor.

\begin{lemma}
	Given $\delta>0$, there exists $\epsilon_{\ref{lem:partition}}=\epsilon_{\ref{lem:partition}}(\delta)>0$ and a positive integer
   $T_0=T_0(\delta)$ such that the following occurs.  Let $T_1,T_2,T_3$ be three positive integers which are divisible by $h$ and with $|T_i-T_j|\in\{0,h\}$, for all
   $i,j\in\{1,2,3\}$ and $T_1>T_0$.  Let $(B_1,B_2,B_3;E)$ be a tripartite graph such that for $\{i,j,k\}=\{1,2,3\}$ $|B_i|=T_j+T_k$, and for $i\ne j$, each vertex in $B_i$ is adjacent to at least $(1-\epsilon_{\ref{lem:partition}})T_1$ vertices in $B_j$.  Then one of the following holds.
\begin{enumerate}
	\item There is a $K_{h,h}$-factor in the graph induced by $(B_1,B_2,B_3;E)$ with the following properties. Each copy is a subgraph of $(B_i,B_j)$ for some $i\ne j$. If we fix a set of at most $\epsilon_{\ref{lem:partition}} T_1$ vertex-disjoint copies of $K_{h,h}$ and at most $\epsilon_{\ref{lem:partition}} T_1$ vertex-disjoint copies of $K_{1, h}$, then  
the $K_{h,h}$-factor contains them as subgraphs.		

	\item The graph induced by $(B_1,B_2,B_3;E)$ can be partitioned such that $B_i=A^{(2)}_i\dotcup A^{(3)}_i$, $|A^{(2)}_i|=T_1$ for $i=1,2,3$ and
   $d(A^{(2)}_j,A^{(2)}_2)\leq\delta$ and
   $d(A^{(3)}_j,A^{(3)}_2)\leq\delta$ for $j=1,3$.
\end{enumerate}
   \label{lem:partition}
\end{lemma}

\begin{figure}
   \begin{center}
      \includegraphics[width=2in]{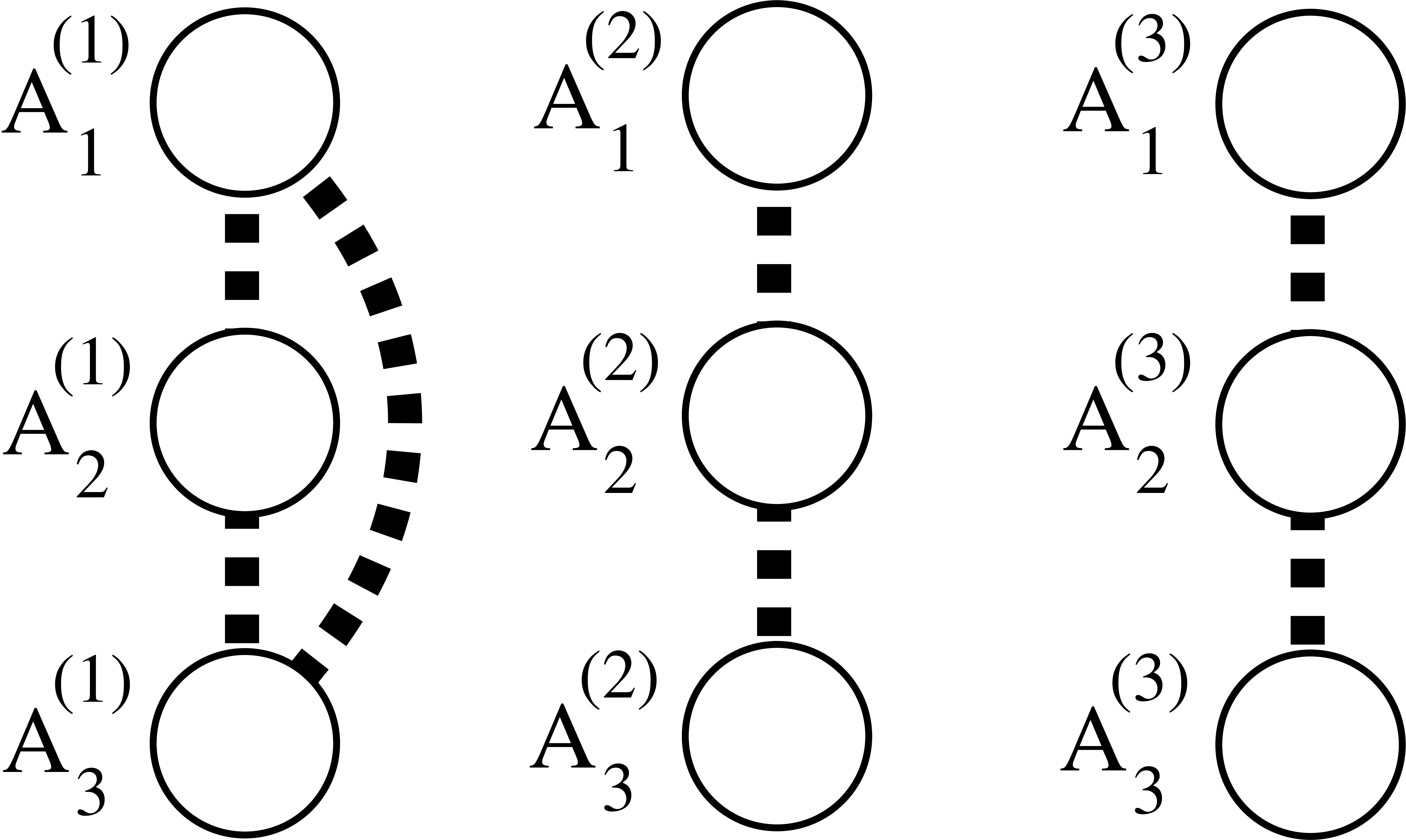}
   \end{center}
   \caption{The diagram that defines Part 2.
   A dotted line represents a sparse pair.}
   \label{fig:figTWO}
\end{figure}

Now to find our $K_{h,h,h}$-factor, we first match vertices in $\widetilde{C}_i$ that are assigned to $B_j$ with $h$ typical neighbors in $B_j$ and those $h$ vertices with $h-1$ typical neighbors in $B_i$.  As the name implies, a typical neighbor is a neighbor which is a typical vertex.  This forms a copy of $K_{h,h,h}$.  Then, place the vertices that were moved in prior steps into copies of $K_{h,h,h}$ by matching the $K_{h,h}$ with vertices in the appropriate ``$\widetilde{A}$'' set.  Remove all of these from $B$, and apply Lemma~\ref{lem:partition} to the remaining adjusted graph with $\delta = \alpha_2$ and $\epsilon_{\ref{lem:partition}} = \alpha_2$.  If the appropriate $K_{h,h}$-factor cannot be found, then we are in the case of Part 2, and $G$ has the form shown in Figure~\ref{fig:figTWO}. A more rigorous definition of this case is provided in Section~\ref{sec:appfigTWO}.

\paragraph{Step~3: Completing the $K_{h,h,h}$-factor.}

If the $K_{h,h}$-factor above is found, then we will recycle notation to define $A^{(1)}_i$ to be the vertices that remain from $\widetilde{A}_i$ after removing copies of $K_{h,h,h}$ as above.  It is easy to see that each $A^{(1)}_i$ will have size close to $T$ and divisible by $h$.  Further define $A^{(j)}_i$, $i=1,2,3$ and $j=2,3$ so that each member of the $K_{h,h}$-factor of $B$ lies in a pair $(A^{(2)}_2,A^{(3)}_3)$, $(A^{(2)}_3,A^{(3)}_1)$ or $(A^{(2)}_1,A^{(3)}_2)$, and so that each of the triples $(A^{(1)}_1,A^{(2)}_2,A^{(3)}_3)$, $(A^{(1)}_2,A^{(2)}_3,A^{(3)}_1)$ and $(A^{(1)}_3,A^{(2)}_1,A^{(3)}_2)$ consist of sets of the same size.  Note that this can be done arbitrarily.

We use Proposition~\ref{prop:sufffactor}, which allows us to
complete a $K_{h,h}$-factor into a $K_{h,h,h}$-factor.  The proof follows easily from K\"onig-Hall.

\begin{proposition} Let $h\geq 1$.
   \begin{enumerate}
      \item \label{it:bi:sufffactor}
      Let $G=(V_1,V_2;E)$ be a bipartite graph with
      $|V_1|=|V_2|=M$, $h$ divides $M$, and each vertex is adjacent to at least $(1-1/(2h^2))M$
      vertices in the other part.  Then, we can find a $K_{h,h}$-factor in $G$.
      \item \label{it:tri:sufffactor}
      Let $G=(V_1,V_2,V_3;E)$ be a tripartite graph with $|V_1|=|V_2|=|V_3|=M$, $h$ divides $M$, and each vertex is adjacent to at least
      $(1-1/(4h^2))M$ vertices in each of the other parts.  Furthermore, let there be a $K_{h,h}$-factor in $(V_2,V_3)$. Then, we can extend it into a $K_{h,h,h}$-factor in $G$.
   \end{enumerate}
   \label{prop:sufffactor}
\end{proposition}

Proposition~\ref{prop:sufffactor}{\it (\ref{it:tri:sufffactor})} allows us to find $K_{h,h,h}$-factors in each of
$(A^{(1)}_1,A^{(2)}_2,A^{(3)}_3)$, $(A^{(1)}_2,A^{(2)}_3,A^{(3)}_1)$ and $(A^{(1)}_3,A^{(2)}_1,A^{(3)}_2)$ which completes the
$K_{h,h,h}$-factor in $G$.

\subsection{Part~2: $G$ is approximately the graph represented by
Figure~\ref{fig:figTWO}} \label{sec:appfigTWO}

Let $H_3$ be the graph on vertices $v_{i}^{(j)}$ for $i,j\in\{1,2,3\}$ with the following non-adjacencies: $v_{1}^{(j)} \not\sim v_{2}^{(j)}$ for $j\in\{1,2,3\}, v_2^{(j)}\not\sim v_3^{(j)}$ for $j\in\{1,2,3\}$ and $v_1^{(1)} \not\sim v_3^{(1)}$.

In this part, our graph $G$ is $\alpha_2$-approximately $H_3$.  It therefore corresponds to the diagram in Figure~\ref{fig:figTWO} in which partite sets are represented as rows, and each row is split into three columns.  Note that the first column of $G$ consists of the pairwise sparse sets from the definition of $\alpha$-extremal, and the second and third columns are defined by the exceptional case of Lemma~\ref{lem:partition}.
We will group the vertices into sets $A_i^{(j)}$ of size between $(1-3\alpha_2^{2/3})T$ and $(1+3\alpha_2^{2/3})T$ so that each vertex in $A^{(j)}_{i}$ is adjacent to at least $\theta T$ vertices in each set $A^{(j')}_{i'}$ when $v_i^{(j)} \sim_{H_3} v_{i'}^{(j')}$.  In other words, the vertices in $A^{(j)}_{i}$ are typical according to the rules established by
Figure~\ref{fig:figTWO}.  The non-typical vertices in row $i$ will be collected in the set $\widetilde{C}_i$.
From this point forward we have issues related to divisibility that we did not have before. Namely, we may need to modify $A^{(2)}_2$ and $A^{(3)}_2$ so that their sizes are divisible by $h$.

\paragraph{Step~1: Ensuring small $A^{(j)}_{i}$ sets of proper size.}

Each $A_{i}^{(j)}$ set has a target size that we will denote $s_{i,j}$.  If $N=3T$, then $s_{i,j}=T$ for all $i,j$.  If $N=3T+h$, let $s_{i,j}=T+h$ for all $i,j$ such that 3 divides $i+j$ and $s_{i,j}=T$ otherwise.  If $N=3T+2h$, let $s_{i,j}=T$ for all $i,j$ such that 3 divides $i+j$ and $s_{i,j}=T+h$ otherwise.  Note that if $N=3T+h$, we can remove one copy of $K_{h,h,h}$ from the triple $(A_1^{(2)}, A_2^{(1)}, A_3^{(3)})$, and if $N=3T+2h$, we can remove two copies of $K_{h,h,h}$ from triples where $i+j$ is not divisible by 3.

Apply Lemma~\ref{lem:stars} to obtain $\max\{0, |A_i^{(j)}|-s_{i,j}\}$ disjoint stars with centers in $A_{i}^{(j)}$ and leaves in $A_{i}^{(j+1)}$.  Then move these star centers to $A_i^{(j')}$ where $j'\ne j$ so that $|A_i^{(j)}|= s_{i,j}$ holds for all $i,j$.




\paragraph{Step~2: Partitioning the sets.}

Before we partition the sets, we must examine the behavior of
$(A^{(2)}_{1}\cup A^{(3)}_{1},A^{(2)}_{3}\cup
A^{(3)}_{3})$.  If this is $\alpha_5$- approximately
$\Theta_{2\times 2}(T)$, then call the \emph{dense} pairs
$(E_1,E_3)$ and $(F_1,F_3)$. Note that the sets $E_i$ and $F_i$ need not be uniquely-defined as long as they satisfy the given condition.  If $(A^{(2)}_{1}\cup A^{(3)}_{1},A^{(2)}_{3}\cup
A^{(3)}_{3})$ is not $\alpha_5$- approximately
$\Theta_{2\times 2}(T)$, do nothing.

For $i\in\{1,3\}$, $\{j,j'\}=\{2,3\}$, we say that $E_i$ and $A_i^{(j)}$ \emph{coincide} if the intersection of their typical vertices is large and therefore the intersection of the typical vertices of $E_i$ and $A_i^{(j')}$ is small.  
We will determine the
quantities that constitute ``large'' and ``small'' later.
If $(E_1,E_3)$
and $(F_1,F_3)$ coincide with
$(A^{(2)}_{1},A^{(3)}_{3})$ and $(A^{(3)}_{1},A^{(2)}_{3})$, respectively, then
$G$ is approximately $\Theta_{3\times 3}(N/3)$.  This case will be handled  in Section~\ref{sec:theta}. If $(E_1,E_3)$ and
$(F_1,F_3)$ coincide with $(A^{(2)}_{1},A^{(2)}_{3})$
and $(A^{(3)}_{1},A^{(3)}_{3})$, respectively, then $G$ is approximately $\Gamma_3(N/3)$.  This case will be handled in Section~\ref{sec:gamma}. Otherwise, there may be no coincidence, or coincidence may occur in exactly one of $V_1$ and $V_3$. Without loss of generality, we will assume that if there is coincidence in only one part, then it occurs in $V_1$.  More specifically, we will assume that $E_1$ coincides with $A_1^{(2)}$, $F_1$ coincides with $A_1^{(3)}$, and neither $E_3$ nor $F_3$ coincides with $A_3^{(j)}$, $j=2,3$.

In addition, note that if, say $A^{(2)}_{1}$
coincides with $E_1$, then every vertex in $A^{(2)}_{1}$ is
adjacent to at least $\theta T$ vertices in $E_3$ and vice
versa. If there is no coincidence, then let $E_1$ and $E_3$ be redefined so that every vertex in $E_1$ is adjacent to at least $\theta T$ vertices in $E_3$ and vice versa. Similarly for $(F_1,F_3)$.

We randomly partition each set $A^{(j)}_{i}$ into two pieces of size divisible by $h$ and as equal as possible.
By the Chernoff bound, with high probability each vertex in $A^{(j)}_{i}$ has at least $(1-2\alpha_3-6\alpha_2^{2/3})(T/2)$
neighbors in each piece of the partition of $A^{(j')}_{i'}$, $i'\neq i$, $j'\neq j$. Moreover, if a vertex has degree at least $\alpha_3 T$ in an $A^{(j)}_{i}$ set, it has degree at least $(\alpha_3/3)(T/2)$ in each of the two partitions.

Let $\Sigma_3$ denote the symmetric group that permutes the elements of $\{1,2,3\}$.  For all $i,j\in\{1,2,3\}$, we assign to each part of $A_i^{(j)}$ a permutation $\sigma\in\Sigma_3$ such that $\sigma(i)=j$ (there are exactly two such permutations) and denote it by $A_{i,\sigma}$.  Furthermore, it is possible to arrange the assignment such that $|A_{1,\sigma}| = |A_{2,\sigma}| = |A_{3,\sigma}|$ for all $\sigma\in\Sigma_3$.
After some adjustments, these permutations will identify which sets the copies of $K_{h,h,h}$ in our covering will span.  For example, a $K_{h,h,h}$ which spans $A_1^{(2)}$, $A_2^{(1)}$ and $A_3^{(3)}$ will be contained in the parts of those sets corresponding to $\sigma=213$, and a $K_{h,h,h}$ which spans $A_1^{(2)}$, $A_2^{(3)}$ and $A_3^{(1)}$ will be contained in the parts of those sets corresponding to $\sigma = 231$.  Note that the permutations $213$ and $231$ are expressed using the notation $\sigma(1)\sigma(2)\sigma(3)$.



\paragraph{Step~3: Assigning vertices.}

Each vertex $c\in \widetilde{C}_2$ has the property that, for all
$j\in\{1,2,3\}$ and distinct $i',i''\in\{1,3\}$, if $c$ is
adjacent to fewer than $\alpha_3 T$ vertices in $A^{(j)}_{i'}$, then $c$ is adjacent to at least $\alpha_3 T$ vertices in $A^{(j)}_{i''}$.

For $i\in\{1,3\}$, each vertex $c\in \widetilde{C}_i$ has the property that, for all $j\in\{1,2,3\}$,  $c$ cannot be adjacent to fewer than $\alpha_3 T$ vertices in either $A^{(2)}_{2}$ or $A^{(3)}_{2}$.  Also, $c$ cannot be adjacent to fewer than $\alpha_3 T$ vertices in both $A^{(1)}_{2}$ and $A^{(1)}_{4-i}$ or both $A^{(2)}_{2}$ and $F_{4-i}$ (if it exists) or both $A^{(3)}_{2}$ and $E_{4-i}$ (if
it exists).  Note that when $F_{4-i}$ and $E_{4-i}$ do not exist, it is because $(A^{(2)}_{1}\cup A^{(3)}_{1},A^{(2)}_{3}\cup
A^{(3)}_{3})$ is not approximately
$\Theta_{2\times 2}(T)$.

Trivially, each vertex in $V_i$ is adjacent to at least
$(1/2-\alpha_3)T$ vertices in at least two of
$\{A^{(1)}_{i'},A^{(2)}_{i'},A^{(3)}_{i'}\}$ and in at least two of $\{A^{(1)}_{i''},A^{(2)}_{i''},A^{(3)}_{i''}\}$, where $i',i''$ are distinct members of $\{1,2,3\}-\{i\}$.  This is particularly important for vertices in $\widetilde{C}_i$.

The $\widetilde{C}_i$ vertices, as well as star-leaves and
star-centers, may only be able to form a $K_{h,h,h}$ with respect
to one particular permutation.

For example, consider a vertex $c$ which had been in $\widetilde{C}_1$ but was put into $A^{(1)}_{1}$ in Step 2.  Then, for either the pair
$(A^{(2)}_{2},A^{(3)}_{3})$ or the pair $(A^{(3)}_{2},A^{(2)}_{3})$, the vertex $c$ is adjacent to at least $(1/2-\alpha_3)T$ vertices in one set and at least $\alpha_3 T$ vertices in the other; otherwise, it would have been a typical vertex in $A^{(1)}_{1}$, $A^{(2)}_{1}$ or $A^{(3)}_{1}$.

Assume that $c$ is adjacent to at least $\alpha_3 T$ vertices in
$A^{(3)}_{2}$ and at least $(1/2-\alpha_3)T$ vertices in $A^{(2)}_{3}$. In this case, if $c$ were placed into the partition corresponding to the identity permutation in Step 3, then exchange $c$ with a vertex in $A_{1,132}$.

In a similar fashion, if there is a star with center in, say
$A^{(2)}_{1}$, and leaves in, say $A^{(1)}_{2}$, then we will use it to form a $K_{h,h,h}$ with respect to the permutation $213\in
\Sigma_3$. Again, if any such leaf or center was placed in the wrong
partition, exchange it with a typical vertex in the other
partition.

The number of leaves in any set is at most
$2h(6\alpha_2^{2/3}T+h)$ and the number of centers is at most $2(6\alpha_2^{2/3}T+h)$; the number of $\widetilde{C}_i$ vertices is at most $9\alpha_2^{2/3}T$. So, if $N$ is large enough, the total number of typical vertices in any $A^{(j)}_{i}$ which were exchanged is at most $2(12h+21)\alpha_2^{2/3}T+4h^2+4h$.

With the partition established and the $\widetilde{C}_i$, star center and leaf vertices in the proper parts, we consider the triple formed by three sets:
\begin{itemize}
   \item $A^{(1)}_{2}$, which will also be denoted $\widetilde{S}_2$
   \item the union of the piece of $A^{(2)}_{1}$ corresponding to $213$ and the piece of $A^{(3)}_{1}$ corresponding to $312$, denoted $\widetilde{S}_1$, and
   \item the union of the piece of $A^{(2)}_{3}$ corresponding to $312$ and the piece of $A^{(3)}_{3}$ corresponding to $213$, denoted $\widetilde{S}_3$.
\end{itemize}
Let the graph induced by the triple
$(\widetilde{S}_1,\widetilde{S}_2,\widetilde{S}_3)$
be denoted $\widetilde{\cal S}$.

\paragraph{Step~4: Finding a $K_{h,h,h}$ cover in $\widetilde{\cal S}$.}

We will first find a $K_{h,h,h}$-factor in $\widetilde{\cal S}$.  This task is complicated because the parts of $\widetilde{\cal S}$ correspond to the permutations $213$ and $312$, meaning the $K_{h,h,h}$'s in our covering either will span $A_1^{(2)}$, $A_2^{(1)}$ and $A_3^{(3)}$ or will span $A_1^{(3)}$, $A_2^{(1)}$ and $A_3^{(2)}$.  If $(A_1^{(2)}\cup A_1^{(3)}, A_3^{(2)}\cup A_3^{(3)})$ is approximately $\Theta_{2\times2}(T)$, then for $i=1,3$, we will need to exchange vertices in $\widetilde{S}_i$ with typical vertices in $E_i$ and $F_i$.  Doing this in the right way will ensure that a $K_{h,h}$-factor of $(\widetilde{S}_1,\widetilde{S}_3)$ can be found, and we will extend that $K_{h,h}$-factor to a $K_{h,h,h}$-factor of $\widetilde{\cal S}$.  Note that this complication does not arise when finding $K_{h,h,h}$'s with respect to permutations in $\Sigma_3 - \{213,312\}$.

To begin, let $T_0=|A^{(1)}_{2}|$.  First, take each existing copy of $K_{1,h}$ in ${\cal \widetilde{S}}$ and complete it to form disjoint copies of $K_{h,h,h}$, using unexchanged typical vertices.  This can be done because $\alpha_4$ is
small enough and the centers are typical vertices. Remove all the copies of $K_{h,h,h}$ that contain stars.

Second, take each vertex $c$ from $\widetilde{C}_i$ and use it to complete a $K_{h,h,h}$.  We can guarantee, because of the random partitioning, that $c$ is adjacent to at least
$(\alpha_3/3)T_0$ vertices in one partition set and $(1/3-2\alpha_3)T_0$
vertices in the other.  Without loss of generality, assume that $c\in
\widetilde{S}_1$ has degree at least $(\alpha_3/3)T_0$ in
$\widetilde{S}_2$ and at least $(1/3-2\alpha_3)T_0$ in
$\widetilde{S}_3$. Since $\alpha_3\gg\alpha_2$, we can guarantee $h$ neighbors of $c$ in $\widetilde{S}_2$ among unexchanged typical vertices and, if $\alpha_3\ll\alpha_4\ll 1$, then $h$ common neighbors of those among unexchanged typical vertices in $N(c)\cap\widetilde{S}_3$. Finally, $\alpha_4\ll h^{-1}$ implies this $K_{h,h}$ has at least $h-1$ more common neighbors in $\widetilde{S}_1$. This is our $K_{h,h,h}$ and we can remove it. Repeat this process for all former members of a $\widetilde{C}_i$.

Third, take each exchanged typical vertex and put it into a
$K_{h,h,h}$ and remove it.  Throughout this process, we have removed at most $C_h\sqrt{\alpha_2}\times T_0$ vertices where $C_h$ is a constant depending only on $h$.  What remains are three sets of the same size, $T'\geq (1-C_h\sqrt{\alpha_2})T_0$, with each vertex in $\widetilde{S}_1$ adjacent to at least, say
$(1/2-2\alpha_4)T'$, vertices in $\widetilde{S}_3$
and vice versa. Each vertex in $\widetilde{S}_1$ and in
$\widetilde{S}_3$ is adjacent to at least
$(1/2-2\alpha_4)T'$ vertices in $\widetilde{S}_2$ and each vertex in $\widetilde{S}_2$ is adjacent to at least
$(1/2-2\alpha_4)T'$ vertices in $\widetilde{S}_1$ and in $\widetilde{S}_3$.

Lemma~\ref{lem:Zhao} (Theorem 9 from \cite{Zhao}) shows that we can find a
factor of $(\widetilde{S}_1,\widetilde{S}_3)$
with vertex-disjoint copies of $K_{h,h}$ unless
$(\widetilde{S}_1,\widetilde{S}_3)$ is
approximately $\Theta_{2\times 2}(T/2)$.
\begin{lemma}[Zhao~\cite{Zhao}]
   For every $\epsilon>0$ and integer $h\geq 1$, there exists
   an $\alpha>0$ and an $N_0$ such that the following holds.
   Suppose that $N>N_0$ is divisible by $h$.  Then every
   bipartite graph $G=(A,B;E)$ with $|A|=|B|=N$ and
   $\delta(G)\geq (1/2-\alpha)N$ either contains a
   $K_{h,h}$-factor, or contains $A'\subseteq A$,
   $B'\subseteq B$ such that $|A'|=|B'|=N/2$ and
   $d(A,B)\leq \epsilon$.
   \label{lem:Zhao}
\end{lemma}

If we can find the factor, apply K\"onig-Hall to form a factor of $\widetilde{\cal S}$ of
vertex-disjoint copies of $K_{h,h,h}$.  If not, apply Lemma~\ref{lem:randpairs}.
Lemma~\ref{lem:randpairs} states, in particular, that if a random
partition results in
$(\widetilde{S}_1,\widetilde{S}_3)$ being
approximately $\Theta_{2\times 2}(T/2)$ with high probability, then
$(A^{(2)}_{1}\cup A^{(3)}_{1},A^{(2)}_{3}\cup
A^{(3)}_{3})$ is approximately $\Theta_{2\times 2}(T)$.  The proof of Lemma~\ref{lem:randpairs} follows from similar arguments to those in the proof of Lemma 3.3 of~\cite{MM} and in Section 3.3.1 of~\cite{MSz} so we omit it.

\begin{lemma}
   For every $\epsilon>0$ and integer $h\geq 1$, there
   exists a $\beta>0$ and positive integer $T_0$ such that
   if $T\geq T_0$ the following holds.  Let $(A,B)$ be a
   bipartite graph such that $|A|,|B|\in\{2T-h,2T,2T+h\}$ with minimum degree at least $(1-\epsilon)T$ and is minimal with respect to this condition.  Let $A'\subset A$, $B'\subset B$,
   $|A'|=|B'|=T$ be chosen uniformly at random.  If
   $$ \Pr\{(A',B')\mbox{ contains a subpair with density
      at most }\epsilon\}\geq 1/4 $$
   then $(A,B)$ is $\beta$-approximately $\Theta_{2\times
   2}(T)$.
   \label{lem:randpairs}
\end{lemma}

We can, therefore, assume the existence of $(E_1,E_3)$ and $(F_1,F_3)$.  Further, we can assume that coincidence occurs only in $V_1$ or not at all; otherwise, we would be in Part 3.

As a result, recall that we let the typical vertices in the dense
pairs in $(A^{(2)}_{1}\cup A^{(3)}_{1},A^{(2)}_{3}\cup
A^{(3)}_{3})$ be denoted $(E_1,E_3)$ and
$(F_1,F_3)$. If the dense pairs do not coincide, then we
will work to ensure that
$|E_1\cap\widetilde{S}_1|=|E_3\cap \widetilde{S}_3|$ and
$|F_1\cap\widetilde{S}_1|=|F_3\cap\widetilde{S}_3|$
and both are divisible by $h$. Do this by moving typical vertices from
$(A^{(2)}_{1}\cap E_1)-\widetilde{S}_1$
into $(A^{(2)}_{1}\cap E_1)\cap\widetilde{S}_1$
and move the same number from $(A^{(2)}_{1}\cap
F_1)\cap\widetilde{S}_1$ into $(A^{(2)}_{1}\cap
F_1)-\widetilde{S}_1$. In addition, move
vertices from $(A^{(2)}_{3}\cap
E_3)-\widetilde{S}_3$ into
$(A^{(2)}_{3}\cap E_3)\cap\widetilde{S}_3$ and
move the same number from $(A^{(2)}_{3}\cap
F_3)\cap\widetilde{S}_3$ into $(A^{(2)}_{3}\cap
F_3)-\widetilde{S}_3$.

This can be done unless one of the intersections $A^{(j)}_{i}\cap
E_i$ or $A^{(j)}_{i}\cap F_i$ is too small.  This implies
the coincidence that we discussed at the beginning of this part. But then, we have guaranteed that the remaining vertices of
$A^{(2)}_{1}$ are not only typical in that set but also typical in $E_1$.  The same is true of $A^{(3)}_{1}$ and $F_1$.

Now, we want to move vertices in $V_3$ to ensure that
$|E_3\cap\widetilde{S}_3|=|A^{(2)}_{1}\cap\widetilde{S}_1|$
and
$|F_3\cap\widetilde{S}_3|=|A^{(3)}_{1}\cap\widetilde{S}_1|$.
Note that we have ensured that both
$|A^{(2)}_{1}\cap\widetilde{S}_1|$ and
$|A^{(3)}_{1}\cap\widetilde{S}_1|$ are divisible by $h$ and
approximately $T/2$.

We can do this as follows: Move vertices from $E_3\cap
A^{(2)}_{3}-\widetilde{S}_3$ to $(E_3\cap
A^{(2)}_{3})\cap\widetilde{S}_3$ and move the same amount from
$(F_3\cap
A^{(2)}_{3})-\widetilde{S}_3$ to $(F_3\cap A^{(2)}_{3})\cap\widetilde{S}_3$.  Also move vertices from
$(E_3\cap A^{(3)}_{3})-\widetilde{S}_3$ to
$(E_3\cap A^{(3)}_{3})\cap\widetilde{S}_3$ and move the same amount from
$(F_3\cap A^{(3)}_{3})-\widetilde{S}_3$ to $(F_3\cap A^{(3)}_{3})\cap\widetilde{S}_3$.  Since none of the intersections are small, this is possible. Moving around these vertices will let us find a $K_{h,h}$-factor of $(\widetilde{S}_1,\widetilde{S}_3)$ which we can complete to a $K_{h,h,h}$-factor of $\widetilde{\cal S}$ by applying Proposition~\ref{prop:sufffactor}{\it (\ref{it:tri:sufffactor})}.

\paragraph{Step~5: Completing the $K_{h,h,h}$-factor in $G$.}

Now that we have found a $K_{h,h,h}$-factor that corresponds to permutations $213$ and $312$, we consider the other permutations in $\Sigma_3$. For a
$\sigma\in\Sigma_3-\{213,312\}$, let ${\cal
S}(\sigma)\stackrel{\rm
def}{=}\left(A_{1,\sigma},A_{2,\sigma},A_{3,\sigma}\right)$
be a triple of parts formed by the random partitioning after the
exchange of vertices has taken place. The set $A_{i,\sigma}$ is a subset of $A^{(\sigma(i))}_{i}$. We have also ensured that $s_{\sigma}\stackrel{\rm
def}{=}\left|A_{1,\sigma}\right|=\left|A_{2,\sigma}\right|=\left|A_{3,\sigma}\right|$ and
$s_{\sigma}$ is divisible by $h$.  It is now easy to ensure that
this triple contains a $K_{h,h,h}$-factor:

First, take each star in ${\cal S}(\sigma)$ and complete it to
form disjoint copies of $K_{h,h,h}$, using unexchanged typical
vertices.  This can be done if $\alpha_4$ is small enough. Remove
all such $K_{h,h,h}$'s containing stars.

Second, take each $c$ which had been a member of some $\widetilde{C}_i$ and
use it to complete a $K_{h,h,h}$.  We can guarantee, because of the
random partitioning, that $c$ is adjacent to at least
$(\alpha_3/3)s_{\sigma}$ vertices in one set and
$(2/3-2\alpha_3)s_{\sigma}$ vertices in the other.  Without loss of
generality, let $c\in A_{1,\sigma}$ with degree at least
$(\alpha_3/3)s_{\sigma}$ in $A_{2,\sigma}$ and at least
$(1/2-2\alpha_3)s_{\sigma}$ in $A_{3,\sigma}$.  Since
$\alpha_3\gg\alpha_2$, we can guarantee $h$ neighbors of $c$ in
$A_{2,\sigma}$ among unexchanged typical vertices and, if
$\alpha_3\ll\alpha_4\ll 1$, then $h$ common neighbors of those among
unexchanged typical vertices in $N(c)\cap A_{3,\sigma}$.
Finally, $\alpha_4\ll h^{-1}$ implies this $K_{h,h}$ has at least
$h-1$ more common neighbors in $A_{1,\sigma}$.  This is our
$K_{h,h,h}$ and we can remove it.  Do this for all former members of
a $\widetilde{C}_i$.

Finally, take each exchanged typical vertex and put it into a
$K_{h,h,h}$ and remove it.  Throughout this process, we have
removed at most $C_h\sqrt{\alpha_2}\times s_{\sigma}$ vertices
where $C_h$ is a constant depending only on $h$.  What remains are
three sets of the same size, $s'\geq
(1-C_h\sqrt{\alpha_2})s_{\sigma}$, with each vertex adjacent to at
least, say $(1-2\alpha_4)s'$, vertices in each of the
other parts. If $N$ is large enough, then we can use the Blow-up
Lemma or
Proposition~\ref{prop:sufffactor}{\it (\ref{it:tri:sufffactor})} to complete the factor of ${\cal S}(\sigma)$ by copies of
$K_{h,h,h}$.

\subsection{Part~3a: $G$ is approximately $\Theta_{3\times
3}(\lfloor N/3\rfloor)$} \label{sec:theta}

Figure~\ref{fig:THETA} shows $\Theta_{3\times 3}$ and we are in the case where $G$ is $\alpha_2$-approximately $\Theta_{3\times 3}(\lfloor N/3\rfloor)$, so $A^{(j)}_i$ and $A^{(j')}_{i'}$ being connected with a dotted line means that the pair $(A^{(j)}_i,A^{(j')}_{i'})$ is sparse.

We will assume for this part that each vertex is adjacent to at
least $h\left\lceil\frac{2N}{3h}\right\rceil+h-1$ vertices in each
of the other pieces of the partition.  Again, let $T=h\lfloor
N/(3h)\rfloor$.

We will group the vertices of $G$ into sets $A_i^{(j)}$ of size between $(1-\sqrt{\alpha_2})T$ and $(1+\sqrt{\alpha_2})T$ so that each vertex in $A_i^{(j)}$ is adjacent to at least $\theta T$
vertices in each set $A^{(j')}_{i'}$ where $i'\neq i$ and $j'\neq j$.  In other words, the vertices in $A_i^{(j)}$ are typical according to the rules established by Figure~\ref{fig:THETA}.  The non-typical vertices in row $i$ will be collected in the set $\widetilde{C}_i$.
Note that each vertex $c\in \widetilde{C}_i$ has the property that, for all
$j\in\{1,2,3\}$ and distinct $i',i''\in\{1,2,3\}-\{i\}$, if
$c$ is adjacent to fewer than $\alpha_3 T$ vertices in
$A^{(j)}_{i'}$, then $c$ is adjacent to at least $\alpha_3 T$
vertices in $A^{(j)}_{i''}$; otherwise $c$ is in some set
$A^{(j)}_{i}$. Furthermore, $c$ is adjacent to at least
$(1/2-\alpha_3)T$ vertices in at least two of
$\{A^{(1)}_{i'},A^{(2)}_{i'},A^{(3)}_{i'}\}$ and in at
least two of
$\{A^{(1)}_{i''},A^{(2)}_{i''},A^{(3)}_{i''}\}$.

\paragraph{Step~1: Ensuring small $A^{(j)}_{i}$ sets of proper size.}

As in Section~\ref{sec:appfigTWO}, each $A_i^{(j)}$ set has a target size $s_{i,j}$.  If $N=3T$, then $s_{i,j}=T$ for all $i,j$.  If $N=3T+h$, let $s_{i,j}=T+h$ when $i=j$ and $s_{i,j}=T$ otherwise.  If $N=3T+2h$, let $s_{i,j}=T$ when $i=j$ and $s_{i,j}=T+h$ otherwise.

Take each triple
$(A^{(j)}_{1},A^{(j)}_{2},A^{(j)}_{3})$, $j=1,2,3$, and
construct disjoint copies of stars so that there are at most $T$
non-center vertices in each set $A^{(j)}_{i}$. We use the fact that
every vertex is adjacent to at least
$h\left\lceil\frac{2N}{3h}\right\rceil+h-1$ vertices in each of the
other parts as well as Lemma~\ref{lem:stars}.  Move these star centers to $A_i^{(j')}$ where $j'\ne j$ so that $|A_i^{(j)}| = s_{i,j}$ holds for all $i,j$.

\paragraph{Step~2: Partitioning the sets.}

We will randomly partition each set $A^{(j)}_{i}$ into two pieces, as
close as possible to equal size but which have size divisible by
$h$, and assign them to a permutation, $\sigma\in \Sigma_3$, which
assigns $\sigma(i)=j$. Each part assigned to
$\sigma$ will be the same size, and these permutations will identify which sets the copies of $K_{h,h,h}$ in our covering
will span.

When $N$ is large, this random partition of $A_i^{(j)}$ will have the following properties with high probability.
A typical vertex in $A^{(j)}_{i}$ has at least
$(1-2\alpha_4-2\sqrt{\alpha_2})(T/2)$ neighbors in each piece of
the partition of $A^{(j')}_{i'}$, $i'\neq i$, $j'\neq j$. Moreover, if a vertex has degree at least
$\alpha_3 T$ in a set, it has degree at least $(\alpha_3/3)(T/2)$
in each of the two partitions.


\paragraph{Step~3: Assigning vertices.}

The $\widetilde{C}_i$ vertices, as well as star centers together with their star-leaves,
may only be able to form a $K_{h,h,h}$ with respect to one
particular permutation.

For example, consider a vertex $c$ which had been in $\widetilde{C}_1$ but
is now in $A^{(1)}_{1}$.  Then, for either the pair
$(A^{(2)}_{2},A^{(3)}_{3})$ or the pair $(A^{(3)}_{2},A^{(2)}_{3})$,
the vertex $c$ is adjacent to at least $(1/2-\alpha_3)T$ in one set
and at least $\alpha_3 T$ vertices in the other.  It is easy to see, since $\alpha_2\ll \alpha_3$, that if this were not true, then $c$ would have been typical with respect to one of the sets $A^{(1)}_{1}$, $A^{(2)}_{1}$ or $A^{(3)}_{1}$, which is a contradiction to the definition of $c$.

Assume that $c$ is adjacent to at least $\alpha_3 T$ vertices in
$A^{(3)}_{2}$ and at least $(1/2-\alpha_3)T$ vertices in $A^{(2)}_{3}$. In
this case, if $c$ were placed into the partition corresponding to
the identity permutation, then exchange $c$ with a typical vertex
in the partition assigned to $132$.

In a similar fashion, if there is a star with center in, say
$A^{(2)}_{1}$, and leaves in, say $A^{(1)}_{2}$, then we will form a
$K_{h,h,h}$ with respect to the permutation $213\in \Sigma_3$.
Again, if any such leaf or center was in the wrong partition,
exchange it with a typical vertex in the other partition.

The number of leaves in any set is at most $2h(\sqrt{\alpha_2}\,
T+h)$ and the number of centers is at most $2(\sqrt{\alpha_2}\,
T+h)$, the number of $\widetilde{C}_i$ vertices is at most $3\sqrt{\alpha_2}\, T$. So, if $N$ is large enough, the total number of typical
vertices in any $A^{(j)}_{i}$ which were exchanged is at most
$(2h+6)\sqrt{\alpha_2}\, T$.

\paragraph{Step~4: Completing the cover.}

For some $\sigma\in \Sigma_3$, let ${\cal S}(\sigma)\stackrel{\rm
def}{=}\left(S^{(\sigma(1))}_{1},S^{(\sigma(2))}_{2},S^{(\sigma(3))}_{3}\right)$
be a triple of parts formed by the random partitioning after the
exchange has taken place. The set $S^{(\sigma(i))}_{i}$ is a subset
of $A^{(\sigma(i))}_{i}$. We have also ensured in Step~3 that
$s_{\sigma}\stackrel{\rm
def}{=}\left|S^{(\sigma(1))}_{1}\right|=\left|S^{(\sigma(2))}_{2}\right|=\left|S^{(\sigma(3))}_{3}\right|$
and $s_{\sigma}$ is divisible by $h$.  It is now easy to ensure that
this triple contains a $K_{h,h,h}$-factor:

First, take each star in ${\cal S}(\sigma)$ and complete it to
form disjoint copies of $K_{h,h,h}$, using unexchanged typical
vertices.  This can be done if $\alpha_4$ is small enough. Remove
all such $K_{h,h,h}$'s containing stars.

Second, take each $c$ which had been a member of some $\widetilde{C}_i$ and
use it to complete a $K_{h,h,h}$.  We can guarantee, because of the
random partitioning, that $c$ is adjacent to at least
$(\alpha_3/3)s_{\sigma}$ vertices in one set and
$(2/3-2\alpha_3)s_{\sigma}$ vertices in the other.  Without loss of
generality, let $c\in S^{(\sigma(1))}_{1}$ have degree at least
$(\alpha_3/3)s_{\sigma}$ in $S^{(\sigma(2))}_{2}$ and at least
$(1/2-2\alpha_3)s_{\sigma}$ in $S^{(\sigma(3))}_{3}$.  Since
$\alpha_3\gg\alpha_2$, we can guarantee $h$ neighbors of $c$ in
$S^{(\sigma(2))}_{2}$ among unexchanged typical vertices and, since
$\alpha_3\ll\alpha_4\ll 1$, $h$ common neighbors of those among
unexchanged typical vertices in $N(c)\cap S^{(\sigma(3))}_{3}$.
Finally, $\alpha_4\ll h^{-1}$ implies this $K_{h,h}$ has at least
$h-1$ more common neighbors in $S^{(\sigma(1))}_{1}$.  This is our
$K_{h,h,h}$ and we can remove it.  Do this for all former members of
a $\widetilde{C}_i$.

Finally, take each exchanged typical vertex and put it into a
$K_{h,h,h}$ and remove it.  Throughout this process, we have
removed at most $\alpha_2^{1/3}s_{\sigma}$ vertices if $\alpha_2$
is small enough. What remains are three sets of the same size,
$s'\geq (1-\alpha_2^{1/3})s_{\sigma}$, with each vertex
adjacent to at least, say $(1-2\alpha_4)s'$, vertices
in each of the other parts. If $N$ is large enough, then we can
use
Proposition~\ref{prop:sufffactor}{\it (\ref{it:tri:sufffactor})} to complete the factor of ${\cal S}(\sigma)$ by copies of
$K_{h,h,h}$.

\subsection{Part~3b: $G$ is approximately $\Gamma_3(\lfloor
N/3\rfloor)$} \label{sec:gamma}

Figure~\ref{fig:GAMMA} shows $\Gamma_3$ and we are in the case where $G$ is $\alpha_2$-approximately $\Gamma_3(\lfloor N/3\rfloor)$, where $A^{(j)}_i$ and $A^{(j')}_{i'}$ being connected with a dotted line means that the pair $(A^{(j)}_i,A^{(j')}_{i'})$ is sparse.

We will assume for this part that each vertex is adjacent to at
least $h\left\lceil\frac{2N}{3h}\right\rceil+h-1$ vertices in each
of the other pieces of the partition.  We also assume that $G$ is
not in \vexc~(see Definition~\ref{def:vexc}).  We must deal with \vexc~separately.

Now, let $T\stackrel{\rm def}{=}h\lfloor N/(3h)\rfloor$.
We may group the vertices of $G$ into sets $A_i^{(j)}$ of size between $(1-\sqrt{\alpha_2})T$ and
$(1+\sqrt{\alpha_2})T$ so that each vertex in $A^{(1)}_{i}$ is adjacent to
at least $(1-\alpha_3)T$ vertices in each set $A^{(j')}_{i'}$ where
$i'\neq i$ and $j'\in\{2,3\}$.
For $j=2,3$, each vertex in $A^{(j)}_{i}$ is adjacent
to at least $(1-\alpha_3)T$ vertices in each set $A^{(1)}_{i'}$ and
$A^{(j)}_{i'}$, where $i'\neq i$.  In other words, the vertices in $A_i^{(j)}$ are typical according to the rules established by Figure~\ref{fig:GAMMA}.  The non-typical vertices in row $i$ will be collected in the set $\widetilde{C}_i$.
Note that each vertex $c\in \widetilde{C}_i$ has the following property: for all
$j\in\{1,2,3\}$ and distinct $i',i''\in\{1,2,3\}-\{i\}$,
if $c$ is adjacent to fewer than $\alpha_3T$ vertices in
$A^{(j)}_{i'}$, then $c$ is adjacent to at least $\alpha_3T$ vertices
in $A^{(j)}_{i''}$. Furthermore, $c$ is adjacent to at least
$(1/2-\alpha_4)T$ vertices in at least two of
$\left\{A^{(1)}_{i'},A^{(2)}_{i'},A^{(3)}_{i'}\right\}$ and
$\left\{A^{(1)}_{i''},A^{(2)}_{i''},A^{(3)}_{i''}\right\}$.

\paragraph{Step~1: Ensuring small $A^{(j)}_{i}$ sets of proper size.}

As in the previous two sections, each $A_i^{(j)}$ has a target size $s_{i,j}$.  There are several cases for $s_{i,j}$ according to the
divisibility of $N/h$. Let $N/h=6q+r$ where $0\leq r<6$.
\begin{itemize}
   \item {\bf $r=0,3$:} $s_{i,j}=T$ for $i=1,2,3$ and $j=1,2,3$.
   \item {\bf $r=1$:} $s_{i,j}=T$ for $i=1,2,3$ and $j=1,3$; and $s_{i,2}=T+h$ for $i=1,2,3$.
   \item {\bf $r=2,5$:} $s_{i,1}=T$ for $i=1,2,3$; and
   $s_{i,j}=T+h$ for $i=1,2,3$ and $j=2,3$.
   \item {\bf $r=4$:} $s_{i,1}=T$ for $i=1,2,3$; and
   $s_{1,3}=s_{2,3}=s_{3,2}=T$; and $s_{1,2}=s_{2,2}=s_{3,3}=T+h$.
\end{itemize}


Without loss of generality, we will assume that both
$|A^{(2)}_{1}|\geq |A^{(3)}_{1}|$ and $|A^{(2)}_{2}|\geq |A^{(3)}_{2}|$.


If $|A^{(2)}_{3}|\geq|A^{(3)}_{3}|$, then $A^{(2)}_{i}$ is larger than
$A^{(3)}_i$ for $i=1,2,3$.  Use
Lemma~\ref{lem:stars}{\it (\ref{lem:stars:bi})} to construct
$\max\left\{|A^{(2)}_{i}|-T,
0\right\}$ disjoint copies of $K_{1,h}$ in the
pair\footnote{Arithmetic in the indices is always done modulo 3.}
$(A^{(2)}_{i},A^{(3)}_{i+1})$ with centers in $A^{(2)}_{i}$. Move these star-centers into $A_i^{(3)}$.

If $|A^{(2)}_{3}|<|A^{(3)}_{3}|$, we do something similar except
that first we use Lemma~\ref{lem:stars}{\it (\ref{lem:stars:bi})} to create the appropriate number of stars in
$(A^{(2)}_{1},A^{(3)}_{2})$ and $(A^{(2)}_{2},A^{(3)}_{1})$ with the
centers in $A^{(2)}_{1}$ and $A^{(2)}_{2}$, respectively. Move these star-centers into $A_1^{(3)}$ and $A_2^{(3)}$, respectively.
Then, after the star-centers have been removed from $A_{2}^{(2)}$, we apply Lemma~\ref{lem:stars}{\it (\ref{lem:stars:bi})} to the pair $(A^{(3)}_{3},A^{(2)}_{2})$, and move the star-centers into $A_3^{(2)}$.

By the conditions on Lemma~\ref{lem:stars}{\it (\ref{lem:stars:bi})}, we see that each remaining set $A^{(j)}_{i}$ is of size at most $T$.  Now, apply Lemma~\ref{lem:stars}{\it (\ref{lem:stars:tri})} to the triple
$(A^{(1)}_{1},A^{(1)}_{2},A^{(1)}_{3})$.  For
star-centers in $A^{(1)}_{i}$, move $T-|A^{(2)}_{i}|$ into $A_i^{(2)}$ 
and $T-|A^{(3)}_{i}|$ into $A_i^{(3)}$. 

If necessary, place vertices from $\widetilde{C}_i$ into
$A^{(j)}_{i}$ for $i=1,2,3$ and $j=1,2,3$, while ensuring that we still
have $|A^{(j)}_{i}|\leq s_{i,j}$.

For $j=2,3$, let ${\cal A}^{(j)}=(A^{(j)}_{1},A^{(j)}_{2},A^{(j)}_{3})$. We remove some copies of $K_{h,h,h}$ from among typical vertices of these sets as follows:
\begin{itemize}
   \item {\bf $r=1$:} One from ${\cal A}^{(2)}$.
   \item {\bf $r=2$:} One from each of ${\cal A}^{(2)}$ and ${\cal
   A}^{(3)}$.
   \item {\bf $r=4$:} One from ${\cal A}^{(2)}$.
   \item {\bf $r=5$:} Two from ${\cal A}^{(2)}$.
\end{itemize}

Recalling $N=(6q+r)h$, each $A_i^{(j)}$ is now of size $2qh$, $2qh+h$ or $2qh+2h$.

\paragraph{Step~2a: Partitioning the sets ($r\neq 3$).}

Let $r\in\{0,1,2,4,5\}$, $\tau_1 = qh$ and $\tau_2 = qh+h$.  Partition each $A^{(j)}_{i}$ set into parts
of nearly equal size.
Each part of the partition will receive a
label $\sigma\in\{1,2,3\}\times\{2,3\}$, where $\sigma = (i,j)$ corresponds to row $i$ and column $j$.  The part with label $(i,j)$ will be denoted $S_{i,j}$. A $K_{h,h,h}$ which is associated with the label $(i,j)$ will span the triple with one part in $A_i^{(1)}$ and two parts in column $j$.  Now, partition each
$A^{(j)}_{i}$ as follows:

Each $A^{(1)}_{i}$ will be split into two pieces.
For $r=0,1,2$ and $i=1,2,3$, both pieces will have size $\tau_1$ and we will
arbitrarily assign the two pieces with the labels $S_{i,2}$ and $S_{i,3}$.
For $r=4$ and $i=3$, assign the piece of size
$\tau_1$ with label $S_{3,3}$ and the one of size $\tau_2$ with $S_{3,2}$.
For $r=4$ and $i=1,2$ and for $r=5$ and $i=1,2,3$, assign the smaller piece with label $S_{i,2}$ and the larger with label $S_{i,3}$.

Each $A^{(2)}_{i}$ will be split into two pieces.  Unless both $r=4$
and $i\in\{1,2\}$, both pieces will be of size $\tau_1$ and will
be assigned $S_{i',2}$ and $S_{i'',2}$ arbitrarily, where
$\{i,i',i''\}=\{1,2,3\}$. If $r=4$ and $i\in\{1,2\}$, the one of
size $\tau_1$ is labeled $S_{3-i,2}$ and the one of size $\tau_2$,
is labeled $S_{3,2}$.

Each $A^{(3)}_{i}$ will be split into two pieces.
If $r\in\{0,1,2\}$, both pieces will be of size $\tau_1$, and if $r=5$ or if $r=4$ and $i=3$, both pieces will be of size $\tau_2$.  In these cases, arbitrarily assign the pieces with labels $S_{i',3}$ and $S_{i'',3}$ where $\{i,i',i''\}=\{1,2,3\}$.
If $r=4$ and $i\in\{1,2\}$, the one of
size $\tau_1$ is labeled $S_{3,3}$ and one of size $\tau_2$ is
labeled $S_{3-i,3}$.

\begin{figure}
\begin{center}
   \epsfig{file=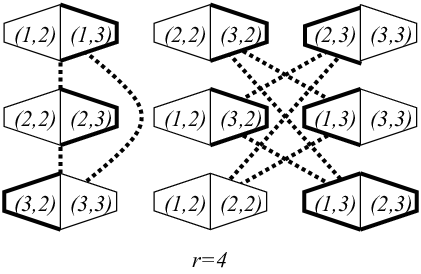}\hspace{.5in}
   \epsfig{file=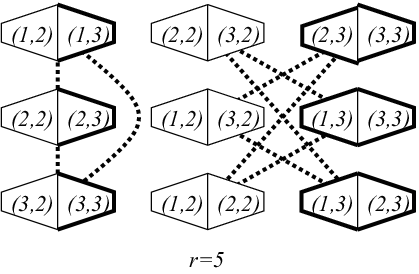}
\end{center}
   \caption{Partitioning the sets.  The light outlined half of a
   set is the piece of size $\tau_1$, the bold outlined half of a
   set is the piece of size $\tau_2$.}
   \label{fig:figGAMtsplit}
\end{figure}

Figure~\ref{fig:figGAMtsplit} diagrams the partitioning for $r=4$ and $r=5$.  Note that when $r\in\{0,1,2\}$, the partition labeling is identical to the case when $r=5$, but all parts have size $\tau_1$.

Partitioning the sets at random again ensures that the above can
be accomplished so that all of the vertices' neighborhoods
maintain roughly the same proportion, as in Part 3a, Step 3.

\paragraph{Step~2b: Partitioning the vertices ($r=3$, not \vexc).}

Let $r=3$ (recall $N=(6q+r)h$) and let $G$ not be in \vexc. We will use Lemma~\ref{lem:stars}{\it (\ref{lem:stars:bi})} to find additional stars between sparse pairs.  Without loss of generality, we seek stars with centers in either $A_1^{(1)}$ or $A_1^{(3)}$.  If we can find at least $h$ centers in one of these sets, then we can make that $A_1^{(j)}$ set of size $2qh$.  If we are not able to do this, every vertex $v\in A_i^{(j)}$ must be adjacent to at most $2h-2$ vertices in $A_{i'}^{(j')}$ where $(A_i^{(j)},A_{i'}^{(j')})$ is a sparse pair.  In turn, we have that every vertex $v\in A_i^{(j)}$ is nonadjacent to at most $3h-3$ vertices in $A_{i'}^{(j')}$ where $(A_i^{(j)},A_{i'}^{(j')})$ is a dense pair.  Since $G$ is approximately $\Gamma_3(\lfloor N/3\rfloor)$, this means $G$ is in \vexc.

Suppose star-centers are removed to make either $|A^{(1)}_{1}|=2qh$ or $|A^{(3)}_{1}|=2qh$.
We will make the set $A^{(2)}_{1}$ of size $(2q+2)h$ by adding star-centers and
vertices from the set $\widetilde{C}_1$.

In each case, if the star-centers that were placed into
$A^{(2)}_{1}$ were themselves originally in $A^{(2)}_{1}$, then we
just treat them as typical vertices again, ignoring the star that
was formed.
Note that all sets are of size $(2q+1)h$, except
$|A^{(2)}_{1}|=(2q+2)h$ and either $A^{(1)}_{1}$ or $A^{(3)}_{1}$,
which has size $2qh$.  If $A^{(1)}_{1}$ is the small set, then
remove one copy of $K_{h,h,h}$ in the triple
$(A^{(3)}_{1},A^{(3)}_{2},A^{(3)}_{3})$.

Now we partition each set as follows: Each $A^{(1)}_{i}$ will have
one piece of size $qh$ with label $(i,3)$. The other set will have
label $(i,2)$ and will be of size $(q+1)h$ in the case of $A^{(1)}_{2}$ and
$A^{(1)}_{3}$ and of size either $qh$ or $(q+1)h$ in the case of
$A^{(1)}_{1}$.  The set $A^{(2)}_{1}$ is partitioned into two pieces
of size $(q+1)h$, one labeled $(2,2)$, the other labeled $(3,2)$.
For $A^{(2)}_{i}$, $i=2,3$, we have one piece of size $qh$ and
labeled $(1,2)$ and the other of size $(q+1)h$, labeled $(5-i,2)$.
For $A^{(3)}_{1}$, it will have two pieces of size $qh$, one labeled
$(2,3)$, the other $(3,3)$. Finally, for $A^{(3)}_{i}$, $i=2,3$, we
have one piece of size $qh$ with label $(5-i,3)$ and the other will
have size either $qh$ or $(q+1)h$ and label $(1,3)$.

Partitioning the sets at random again ensures that the above can
be accomplished so that all of the vertices' neighborhoods
maintain roughly the same proportion, as in Part 3a, Step 3.

\paragraph{Step~3: Assigning vertices.}

For any $\sigma\in\{1,2,3\}\times\{2,3\}$, we will show that the
star-centers and $\widetilde{C}_i$ vertices, in any $A^{(j)}_{i}$ can be assigned to one
of the two parts of the partition.

For example, consider a vertex $c$ which had been in $\widetilde{C}_1$ but
is now in $A^{(1)}_{1}$.  Then, for either the pair
$(A^{(2)}_{2},A^{(2)}_{3})$ or the pair $(A^{(3)}_{2},A^{(3)}_{3})$,
the vertex $c$ is adjacent to at least $(1/2-\delta)T$ in one set
and at least $\alpha_3 T$ vertices in the other.  If such a pair is
$(A^{(2)}_{2},A^{(2)}_{3})$ then if $c$ were labeled $(1,2)$
exchange it with a typical vertex with label $(1,3)$.

Now, for example, consider a vertex $c$ which had been in $\widetilde{C}_1$
but is now in $A^{(2)}_{1}$.  It is easy to check that for either
the pair $(A^{(1)}_{2},A^{(2)}_{3})$ or the pair
$(A^{(1)}_{3},A^{(2)}_{2})$, the vertex $c$ is adjacent to at least
$(1/2-\alpha_3)T$ in one set and at least $\alpha_3 T$ vertices in
the other.  If such a pair is, say, $(A^{(1)}_{2},A^{(2)}_{3})$, and
$c$ is not labeled $(2,2)$, then exchange it for a typical vertex of
that label.

A similar analysis can be applied to any $c\in \widetilde{C}_i$ for $i=1,2,3$.

Now we consider stars.  All star-centers are in sets $A^{(2)}_{i}$ or
$A^{(3)}_{i}$.  Without loss of generality, assume $z$ is such a
center in $A^{(2)}_{1}$ and the leaves are in $V_2$.  If the leaves
are in $A^{(1)}_{2}$, then $z$ must have been a member of $A^{(1)}_{1}$
originally.  So, $z$ and its leaves must have label $(2,2)$.  If
the leaves are in $A^{(2)}_{2}$, then $z$ must have been a member of
$A^{(3)}_{1}$ originally.  So, $z$ and its leaves must have label
$(3,2)$.  Exchange $z$ with typical vertices to ensure this.

Finally, we consider typical vertices moved from $A^{(2)}_{i}\cup
A^{(3)}_{i}$ to $A^{(1)}_{i}$.  Without loss of generality, suppose $z$ is
such a vertex in $A^{(1)}_{1}$.  If $z$ were originally from
$A^{(2)}_{1}$, then it is a typical vertex with respect to $A^{(2)}_{2}$
and $A^{(2)}_{3}$ and $z$ should receive label $(1,2)$.  Otherwise, it
is typical with respect to $A^{(3)}_{2}$ and $A^{(3)}_{3}$ and $z$ should
receive label $(1,3)$.

This completes the verification that all moved vertices can
receive at least one label of the $A^{(j)}_{i}$ set in which it is
placed.

\paragraph{Step~4: Completing the cover.}

For any $\sigma\in\{1,2,3\}\times\{2,3\}$, let ${\cal S}(\sigma)$
be the triple of parts with label $\sigma$.  Note that the label $(i,j)$ corresponds to a triple with one part in $A_i^{(1)}$ and two parts in column $j$.  We can finish the $K_{h,h,h}$-factor as in Part 3a,
Step 5.


\subsection{\Vexc}
\label{sec:vexc}

Recall \vexc:
\begin{quote}
   There are integers $N,q$ such that $N=(6q+3)h$.  There are
   sets $A^{(j)}_{i}$ for $i,j\in\{1,2,3\}$, with sizes at least
   $2qh+1$, such that if $v\in A^{(j)}_{i}$ then $v$ is nonadjacent
   to at most $3h-3$ vertices in $A^{(j')}_{i'}$ whenever the pair
   $(A^{(j)}_{i},A^{(j')}_{i'})$ corresponds to an edge in the graph
   $\Gamma_3$ with respect to the usual correspondence.
\end{quote}

In this case, we must raise the minimum degree condition to
$2N/3+2h-1$.  Recalling Part 4, Step 3b, we were able to proceed if
we were able to make one of the sets $A^{(j)}_{i}$ small by means of
creating stars.  Each vertex in $A^{(2)}_{2}$ is adjacent to at
least $|A_1^{(3)}|-N/3+2h-1$ vertices in $A^{(3)}_{1}$.  Using
Lemma~\ref{lem:stars}{\it (\ref{lem:stars:bi})}, we have that there is a
family of $|A^{(3)}_{1}|-N/3+h$ vertex-disjoint stars with centers
in $A^{(3)}_{1}$.  We move the centers to $A^{(2)}_{1}$.  Then we
can proceed from Part 3b, Step 4.


\subsection{Proofs of Lemmas}
\label{sec:lemmas}

Lemma~\ref{lem:stars} is used to find vertex-disjoint $h$-stars in a graph $G$.  Part~{\it (\ref{lem:stars:bi})} deals with the case where $G=({A}_1,{A}_2;E)$, and part~{\it (\ref{lem:stars:tri})} deals with the case where $G=({A}_1,{A}_2,{A}_3;E)$. \\

\begin{proofcite}{Lemma~\ref{lem:stars}}
   \begin{enumerate}[label=(\arabic*),font=\itshape]
      \item Let $\delta_1=d_1-h+1$.  If the stars cannot be
      created greedily, then
      there is a set $S\subset {A}_1$ and a set $T\subset {A}_2$ such
      that $|S|\leq\delta_1-1$ and $|T|=|S|h$ and each
      vertex in ${A}_1- S$ is adjacent to at most $h-1$
      vertices in ${A}_2- T$.  In this case,
      $$ (d_1-|S|)|{A}_2- T|\leq e({A}_1- S,
      {A}_2- T)\leq (h-1)|{A}_1- S| . $$

      This gives
      \begin{eqnarray*}
         |S| & \geq & \delta_1
                      -(h-1)\frac{|A_1- S|-|A_2- T|}
                                 {|A_2- T|} \\
               & = & \delta_1-(h-1)\frac{|A_1|-|A_2|+(h-1)|S|}
                                           {|A_2|-h|S|} \\
               & \geq & \delta_1-(h-1)\frac{(h+1)\epsilon M}
                                           {(1-(h+1)\epsilon)M} .
      \end{eqnarray*}
      If $\epsilon<(h^2+h)^{-1}$, then this gives
      $|S|>\delta_1-1$.  Since $|S|$ is an integer,
      $|S|\geq\delta_1$, contradicting the condition we put on
      $|S|$.

      \item Let $\delta_i=\max\{0,d_i-h+1\}$ for $i=1,2,3$.  If,
      say, $\delta_3=0$, then apply part {\it (\ref{lem:stars:bi})} to
      the pair $(A_2,A_3)$ to create $\delta_2$ vertex-disjoint
      stars with centers in $A_2$.  Let $Z_2$ be the set of the
      centers.  Apply part {\it (\ref{lem:stars:bi})} to
      $(A_1,A_2- Z_2)$ and we can find $\delta_1$
      vertex-disjoint stars with centers in $A_1$ if
      $2\epsilon<(h^2+h)^{-1}$.

      So, we may assume that $\delta_i>0$ for $i=1,2,3$.  Note that
      if it is possible to construct $\delta_1+\delta_2$ disjoint
      copies of $K_{1,h}$ in $(A_1,A_2)$ with centers,
      $Z_1\subset A_1$, then we can finish by applying part
      {\it (\ref{lem:stars:bi})}.  To see this, apply part
      {\it (\ref{lem:stars:bi})} to $(A_3,A_1- Z_1)$,
      with
      $3\epsilon<(h^2+h)^{-1}$, creating $\delta_3$ stars with
      centers $Z_3\subset A_3$.  Then apply part {\it (\ref{lem:stars:bi})}
      to $(A_2,A_3- Z_3)$. (Here, we need $2\epsilon<(h^2+h)^{-1}$.)
      There will be $\delta_1$ stars remaining in $(A_1,A_2)$ which
      are vertex-disjoint from the rest.

      So, we will assume that it is not possible to create
      $\delta_1+\delta_2$ vertex-disjoint copies of $K_{1,h}$
      in $(A_1,A_2)$ with centers in $A_1$.  That means there is an
      $S\subset A_1$ and a $T\subset A_2$ such that
      $|S|<\delta_1+\delta_2$, $|T|=h|S|$ and every vertex in
      $A_1- S$ is adjacent to at most $h-1$ vertices in
      $A_2- T$.

      Now apply part {\it (\ref{lem:stars:bi})} to $(A_3,A_1- S)$
      to obtain $\delta_3$ vertex-disjoint copies of $K_{1,h}$ with
      centers $Z_3\subset A_3$.  (Here, we need
      $3\epsilon<(h^2+h)^{-1}$.) Next, apply part
      {\it (\ref{lem:stars:bi})} to $(A_2,A_3- Z_3)$ to obtain
      $\delta_2$ vertex-disjoint copies of $K_{1,h}$ with centers
      $Z_2 \subset A_2$.  (Here, we need
      $2\epsilon<(h^2+h)^{-1}$.)  Finally, apply part
      {\it (\ref{lem:stars:bi})} to $\left(A_1,A_2- (Z_2\cup
      T)\right)$ to obtain $\delta_1$ vertex-disjoint copies of
      $K_{1,h}$ with centers $Z_1\subset A_1$.  (Here, we need
      $(2h+2)\epsilon<(h^2+h)^{-1}$.)  But, because no vertex in
      $A_1- S$ is adjacent to $h$ vertices in
      $A_2- (Z_2\cup T)$, it must be the case that
      $Z_1\subset S$ and our $\delta_1+\delta_2+\delta_3$ copies of
      $K_{1,h}$ are, indeed, vertex-disjoint.
   \end{enumerate}
\end{proofcite}

Lemma~\ref{lem:superstars} is used to find a copy of $K_{1,h,h}$ in a tripartite graph $(B_1,B_2,B_2;E)$.  If a $K_{1,h,h}$ cannot be found, then the graph must be approximately $\Theta_{3\times 2}(M)$. \\

\begin{proofcite}{Lemma~\ref{lem:superstars}}
   We can first apply the following theorem of
   Erd\H{o}s, Frankl and R\"odl~\cite{EFR}:
   \begin{theorem}
      For every $\epsilon'>0$ and graph $F$, there is a constant
      $n_0$ such that for any graph $G$ of order $n\geq n_0$, if
      $G$ does not contain $F$ as a subgraph, then $G$ contains a
      set $E'$ of at most $\epsilon' n^2$ edges such that
      $G- E'$ contains no $K_r$ with $r=\chi(F)$.
      \label{thm:erdos}
   \end{theorem}
   Here, $F=K_{1,h,h}$ and $r=3$.

   Let us remove at most $\epsilon' (6M)^2$ edges from $G$ so that it becomes triangle-free.  In doing so, some vertices might be nonadjacent to many more vertices than before.  We want to remove such vertices so that we can apply Proposition~\ref{pTHETA}, which appeared in~\cite{MM} and is rephrased below:
	\begin{proposition}
      For a $\Delta$ small enough, there exists $\epsilon''>0$
      such that if $H$ is a tripartite graph with at least
      $2\left(1-\epsilon''\right)t$ vertices in each vertex class
      and each vertex is nonadjacent to at most
      $\left(1+\epsilon''\right)t$ vertices in each of the other
      classes. Furthermore, let $H$ contain no triangles.  Then,
      each vertex class is of size at most
      $2\left(1+\epsilon''\right)t$ and $H$ is $\Delta$-approximately
      $\Theta_{3\times 2}(t)$. \label{pTHETA}
   \end{proposition}

   For $\epsilon''\gg\epsilon'$, at least $2(1-\epsilon'')M$ vertices are nonadjacent to at most $(1+\epsilon'')M$ vertices in each of the other classes.  Otherwise, we would have had to delete a total of at least $\Omega(\epsilon'')M$ edges incident to each of these vertices, of which there would be at least $\Omega(\epsilon'')M$.  But this means deleting $\Omega((\epsilon'')^2)M^2$ edges, which is more than $\epsilon'(6M)^2$.

   So we apply Proposition~\ref{pTHETA}.  Thus, $G$ is approximately $\Theta_{3\times 2}(M)$, and so the lemma follows.
\end{proofcite}

Lemma~\ref{lem:partition} is used to find a $K_{h,h}$-factor in a tripartite graph $(B_1,B_2,B_2;E)$.  If the factor cannot be found, then the graph has a structure like columns $2$ and $3$ of the diagram in Figure~\ref{fig:figTWO}.  

~\\

\begin{proofcite}{Lemma~\ref{lem:partition}}
	Let $\epsilon'$ be chosen such that $\epsilon'\ll\delta$.

   For this lemma, we partition the possibilities according to
   whether the pairs $(B_i,B_j)$ are approximately
   $\Theta_{2\times 2}(T_1)$.  That is, there are two pairs of sets of
   size $T_1$ which have density less than $\epsilon'$.
   Minimality gives the rest.

   In addition, we say that graphs $\Theta_{2\times 2}(T_1)$
   \textdef{coincide} if $(B_i,B_j)$ and  $(B_j,B_k)$ are approximately  $\overline{B}_i\subseteq B_i$,
   $\overline{B}_j\subseteq B_j$, $\overline{B}_k\subseteq B_k$, all of size
   $T_1$, such that both $(\overline{B}_i,\overline{B}_j)$ and $(\overline{B}_j,\overline{B}_k)$ have density less than $\epsilon'$.  Note that this means that $(B_i- \overline{B}_i,B_j- \overline{B}_j)$ and $(B_i- \overline{B}_i,B_j- \overline{B}_j)$

   \paragraph{Case~1: No pair is $\Theta_{2\times 2}(T_1)$.}

   For each distinct $i,j,k\in\{1,2,3\}$, partition $B_i$ into   two pieces, $B_i[j]$ and $B_i[k]$ with $|B_i[j]|=T_j$ and $|B_i[k]|=T_k$.  If this partition is done uniformly at random, then with probability approaching 1, each vertex in $B_i[k]$ is adjacent to at least $(1/2-\epsilon^{1/2})T_k$
   vertices in $B_j[k]$.  So there exists a partition such
   that each vertex in $B_i$ is adjacent to at least
   $(1/2-\epsilon^{1/2})T_1$ vertices in each of the pieces
   $B_j[k]$, $j,k\neq i$ and such that the pair
   $(B_2[1],B_3[1])$ fails to contain a subpair with
   $\lfloor T_1/2\rfloor$ vertices in each part and density at
   most $\epsilon^{1/3}$.

   The vertices that are reserved will have to be placed in the
   proper set.  For example, if a reserved $K_{h,h}$ is in the
   pair $(B_i,B_j)$, then those vertices will need to be in the
   pair $(B_i[k],B_j[k])$.  So, we exchange vertices in
   $B_i[k]$ for vertices in $B_i[j]$ so that reserved
   vertices are in the proper place.  At most
   $4(\epsilon+\epsilon)T_1$ vertices are either reserved or
   moved in each set $B_i[j]$.  After such exchanges occur,
   place the moved vertices into vertex-disjoint copies of $K_{h,h}$ that lie entirely within the given pairs.  This can be done because each vertex not in $B_i$ is adjacent to almost half of the vertices in both $B_i[j]$ and $B_i[k]$.

   Consider what remains of these sets.  The number of vertices
   is still divisible by $h$ and at most $8h(\epsilon)T_1$ have been placed into these copies of $K_{h,h}$.  We look for a perfect $K_{h,h}$-factor in each of the pairs $(B_1[3],B_2[3])$, $(B_1[2],B_3[2])$ and $(B_2[1],B_3[1])$.  Recall that each of these pairs has minimum degree at least $(1/2-\epsilon^{1/2})T_1$.  Utilizing a lemma in~\cite{Zhao} --  stated as Lemma~\ref{lem:Zhao} in Section~\ref{sec:appfigTWO} above -- we are able to find such a factor unless at least one of those pairs is $\alpha(\epsilon^{1/2})$-approximately
   $\Theta_{2\times 2}(T_1/2)$. (Minimality gives the other
   sparse pair.)

   Lemma~\ref{lem:randpairs} says that if random selections give a graph that is approximately $\Theta_{2\times 2}(T_1/2)$, then the original graph was, too.  So, along with Lemma~\ref{lem:Zhao}, it establishes that if, after moving our vertices, we are unable to complete our $K_{h,h}$-cover in $(B_i[k],B_j[k])$ with nontrivial probability, then the pair $(B_i,B_j)$ is $\epsilon'$-approximately $\Theta_{2\times 2}(T_1)$, where $\epsilon'=\beta(\alpha(\epsilon^{1/2}))$.

   Since none of the pairs is $\epsilon'$-approximately
   $\Theta_{2\times 2}(T_1)$, we can find the required factor
   of $\left(B_1,B_2,B_3\right)$ by copies of $K_{h,h}$.

\paragraph{Case~2: Exactly one pair is $\Theta_{2\times 2}(T_1)$.}

   Here, we will assume that $B_1=\overline{B}_1\dotcup \widehat{B}_1$ and
   $B_2=\overline{B}_2\dotcup \widehat{B}_2$, where
   $|\overline{B}_1|=|\widehat{B}_2|=T_1$ and
   $d(\overline{B}_1,\widehat{B}_2), d(\widehat{B}_1,\overline{B}_2)\leq\epsilon'$.
   A random partition of $B_1$ into pieces, with probability approaching 1 as $T_1$ approaches infinity, will partition $\overline{B}_1$ into two approximately equal pieces.  In particular, let the \textdef{typical vertices} in $\overline{B}_1$ be those that are nonadjacent to at
   most $(\epsilon')^{1/2}T_1$ in $\widehat{B}_2$.  There are at most $(\epsilon')^{1/2}T_1$ such vertices.  A similar conclusion can be drawn from $\overline{B}_2$, $\widehat{B}_1$ and $\widehat{B}_2$.

   In this case, we randomly partition $B_1$, $B_2$ and $B_3$ into the sets $B_i[k]$ as prescribed.  Exchange the vertices as we have done above and complete both the reserved and exchanged vertices to form copies of $K_{h,h}$.  This encompasses at most $8h\epsilon T_1$ vertices.
   Exchange vertices in $B_1[3]$ with vertices in $B_1[2]$ and vertices in $B_2[3]$ with vertices in $B_2[1]$ so that there are exactly $h\lfloor T_1/(2h)\rfloor$ typical vertices of $\overline{B}_1$ in $B_1[3]$ and $h\lfloor T_1/(2h)\rfloor$ typical vertices of $\widehat{B}_2$ in $B_2[3]$.  Let the rest of the vertices, not matched into a $K_{h,h}$, in $B_1[3]$ be typical vertices in $\widehat{B}_1$ and the rest of the vertices in $B_2[3]$ be typical in $\overline{B}_2$.  Using
   Proposition~\ref{prop:sufffactor}{\it (\ref{it:bi:sufffactor})} on
   each pair of sets of typical vertices in $(B_1[3],B_2[3])$ will easily have a $K_{h,h}$-factor.  With $\epsilon'$ small enough, we can guarantee that at most $(\epsilon')^{1/3}T_1$ vertices in $(B_1[2],B_3[2])$ and $(B_2[1],B_3[1])$ were moved.  Applying
   Lemmas~\ref{lem:Zhao} and~\ref{lem:randpairs}, and the fact that no pair other than $(B_1,B_2)$ can be $\epsilon'$-approximately
   $\Theta_{2\times 2}(T_1)$, we conclude that the pairs
   $(B_1[2],B_3[2])$ and $(B_2[1],B_3[1])$
   can be completed to $K_{h,h}$-factors.

\paragraph{Case~3: Exactly two pairs are $\Theta_{2\times 2}(T_1)$, which do not coincide.}

   Let the pairs in question be $(B_1,B_2)$ and $(B_2,B_3)$. Let the dense pairs in the subgraph induced by $(B_1,B_2)$ be $(\overline{B}_1,\overline{B}_2)$ and
   $(\widehat{B}_1,\widehat{B}_2)$.  Let the dense pairs in $(B_2,B_3)$ be $(\mathring{B}_2,\mathring{B}_3)$ and
   $(\ddot{B}_2,\ddot{B}_3)$.  Moreover, since the pairs
   fail to coincide, we can conclude that the intersection of the typical vertices of $\overline{B}_2$ with the typical vertices of each of $\mathring{B}_2$ and $\ddot{B}_2$ is at least $(\epsilon')^{1/4}T_1$ and similarly for $\widehat{B}_2$.

   Once again, we randomly partition the vertices in $B_1$, $B_2$ and $B_3$ and move vertices so as to ensure that the reserved vertices and the vertices exchanged for them are placed into vertex-disjoint copies of $K_{h,h}$.  Our concern at this point is the vertices in $B_2$.

   Consider the vertices in $(B_1[3],B_2[3])$.
   Approximately half are typical vertices of $\overline{B}_2$ and approximately half are typical vertices of $\widehat{B}_2$.  Take each non-typical vertex in $B_1[3]$ and in $B_2[3]$, match them with a copy of $K_{h,h}$ in the pair $(B_1[3],B_2[3])$ and remove them.  Do the same for vertices in $B_2[1]$ that are not typical in $\mathring{B}_2$ or $\ddot{B}_2$ and in $B_3[1]$ that are not typical in $\mathring{B}_3$ or $\ddot{B}_3$. Remove those copies of $K_{h,h}$ also.

   Observe that there are at least $\epsilon^{1/4}t_1/4$ vertices
   in each intersection of $\overline{B}_2$ or $\widehat{B}_2$ with $\mathring{B}_2$ or $\ddot{B}_2$ and with $B_2[3]$ or $B_2[1]$.

   First, move $a$ vertices from $\overline{B}_2\cap
   \mathring{B}_2\cap B_2[3]$ to $\overline{B}_2\cap
   \mathring{B}_2\cap B_2[1]$ to make $|\overline{B}_2\cap B_2[3]|$ divisible by $h$.  Second, move $a+b$ vertices from
   $\widehat{B}_2\cap \mathring{B}_2\cap B_2[1]$ to
   $\widehat{B}_2\cap \mathring{B}_2\cap B_2[3]$ to
   make $|\mathring{B}_2\cap B_2[1]|$ divisible by $h$.  Third, move $a+b+c$ vertices from $\widehat{B}_2\cap \ddot{B}_2\cap B_2[3]$ to $\widehat{B}_2\cap \ddot{B}_2\cap B_2[1]$. This will make both $|\widehat{B}_2\cap B_2[3]|$ and
   $|\ddot{B}_2\cap B_2[1]|$ divisible by $h$.

   Here $a$, $b$ and $c$ are the remainders of $|\overline{B}_2\cap B_2[3]|$, $|\mathring{B}_2\cap B_2[1]|$ and $|\widehat{B}_2\cap B_2[3]|$, respectively, when each is divided by $h$.  Observe that both $|\overline{B}_2\cap B_2[3]|+|\widehat{B}_2\cap B_2[3]|$ and $|\mathring{B}_2\cap B_2[1]|+|\ddot{B}_2\cap B_2[1]|$ are divisible by $h$.

   Finally, we exchange vertices in $\overline{B}_1\cap B_1[3]$ with those in $\overline{B}_1\cap B_1[2]$ so that $|\overline{B}_1\cap B_1[3]|=|\overline{B}_2\cap B_2[3]|$ and similarly for $\widehat{B}_2$.  Also, exchange vertices
   in $\mathring{B}_3\cap B_3[1]$ with those in
   $\mathring{B}_3\cap B_3[2]$ so that
   $|\mathring{B}_3\cap B_3[1]|=|\mathring{B}_2\cap
   B_2[1]|$ and similarly for $\ddot{B}_2$.

   Then, in $(\overline{B}_1\cap B_1[3],\overline{B}_2\cap
   B_2[3])$, first greedily place each moved vertex into copies of $K_{h,h}$ and then finish the factor via
   Proposition~\ref{prop:sufffactor}{\it (\ref{it:bi:sufffactor})}.
   Do the same for $\left(\widehat{B}_1\cap B_1[3],\widehat{B}_2\cap B_2[3]\right)$,
   $\left(\mathring{B}_2\cap B_2[1],\mathring{B}_3\cap B_3[1]\right)$ and $\left(\mathring{B}_2\cap
   B_2[1],\mathring{B}_3\cap B_3[1]\right)$.

   Finally, we can complete the factor of $(B_1[2],B_3[2])$ because if it is not possible, Lemmas~\ref{lem:Zhao} and~\ref{lem:randpairs} would require $(B_1,B_3)$ to be approximately $\Theta_{2\times 2}(T_1)$, excluded by this case.

\paragraph{Case~4: Three pairs are $\Theta_{2\times 2}(T_1)$, none of which coincide.}

   Let the dense pairs in $(B_1,B_2)$ be
   $(\overline{B}_1,\overline{B}_2)$ and
   $(\widehat{B}_1,\widehat{B}_2)$.
   Let the dense pairs in $(B_2,B_3)$ be
   $(\mathring{B}_2,\mathring{B}_3)$ and
   $(\ddot{B}_2,\ddot{B}_3)$.  Let the dense pairs in $(B_1,B_3)$ be $(B_1^{\sharp},B_3^{\sharp})$ and $(B_1^{\flat},B_3^{\flat})$.  Moreover,
   since the pairs fail to coincide, we can conclude that the
   intersection of the typical vertices of one set of sparse pairs with the typical vertices of another is at least
   $(\epsilon')^{1/4}T_1$.

   Partition $B_1$, $B_2$ and $B_3$ into
   appropriately-sized sets as before, uniformly at random.  The degree conditions hold with high probability as before.  Take non-typical vertices and complete them greedily to place them in vertex-disjoint copies of $K_{h,h}$ within each of the pairs $(B_1[3],B_2[3])$, $(B_2[1],B_3[1])$ and $(B_1[2],B_3[2])$.  Remove these copies of $K_{h,h}$ from the graph.

   Let $M$ be the largest multiple of $h$ less than or equal to
   the size of the intersection of what remains of any sparse set
   (\ie, $\overline{B}_i,\widehat{B}_i,
   \mathring{B}_i,\ddot{B}_i, B_i^{\sharp},B_i^{\flat}$) with a set of the form
   $B_i[k]$.

   We can move vertices as in Case 3 by letting
   $a=|\overline{B}_2\cap B_2[3]|-M$,
   $b=|\mathring{B}_2\cap B_2[1]|-M$ and
   $c=|\widehat{B}_2\cap B_2[3]|+M-T_3$, which is also equal to $T_1-M-a-b-|\ddot{B}_2\cap B_2[1]|$.  We can
   perform similar operations to guarantee that, among the vertices that remain in the graph, that
   \begin{eqnarray*}
   M & = & \left|\overline{B}_1\cap
   B_1[3]\right|=\left|\overline{B}_2\cap B_2[3]\right|=\left|\mathring{B}_2\cap
   B_2[1]\right|=\left|\mathring{B}_3\cap
   B_3[1]\right| \\
   & = & \left|B_1^{\sharp}\cap
   B_1[2]\right|=\left|B_3^{\sharp}\cap B_3[2]\right|
   \end{eqnarray*}
   The fact that the pairs do not coincide ensures that there are
   enough vertices to make these moves.

   Place the moved vertices into vertex-disjoint copies of
   $K_{h,h}$ and finish the factor via
   Proposition~\ref{prop:sufffactor}{\it (\ref{it:bi:sufffactor})}.

\paragraph{Case~5: There are at least two pairs which are $\Theta_{2\times 2}(T_1)$ and which coincide.}

   This is exactly the exceptional case stated in the lemma and without loss of generality the pairs $(A^{(2)}_1,A^{(2)}_2)$ and $(A^{(2)}_2,A^{(2)}_3)$ are those that witness the coincidence of the copies of $\Theta_{2\times 2}(T_1)$.

\end{proofcite}




\begin{acknowledgements}
The authors would like to acknowledge and thank the Department of Mathematics, Statistics, and Computer Science at the University of Illinois at Chicago for their supporting Martin via a visitor fund.  The authors also wish to acknowledge the support of National Security Agency grant H98230-08-1-0015 for Hogenson for a summer research assistantship.
\end{acknowledgements}


\bibliographystyle{spmpsci}      
\bibliography{3tilebib}   

\end{document}